\pdfoutput=1
\RequirePackage{ifpdf}
\ifpdf 
\documentclass[pdftex]{sigma}
\else
\documentclass{sigma}
\fi

\numberwithin{equation}{section}

\newtheorem{Theorem}{Theorem}[section]
\newtheorem{Corollary}[Theorem]{Corollary}
\newtheorem{Lemma}[Theorem]{Lemma}
\newtheorem{Proposition}[Theorem]{Proposition}
 { \theoremstyle{definition}
\newtheorem{Definition}[Theorem]{Definition}
\newtheorem{Note}[Theorem]{Note}
\newtheorem{Example}[Theorem]{Example}
\newtheorem{Remark}[Theorem]{Remark} }

\begin{document}

\allowdisplaybreaks

\renewcommand{\PaperNumber}{047}

\FirstPageHeading

\ShortArticleName{The Universal Askey--Wilson Algebra and DAHA of Type $(C_1^{\vee},C_1)$}

\ArticleName{The Universal Askey--Wilson Algebra
\\
and DAHA of Type $\boldsymbol{(C_1^{\vee},C_1)}$}

\Author{Paul TERWILLIGER}

\AuthorNameForHeading{P.~Terwilliger}

\Address{Department of Mathematics, University of Wisconsin, Madison, WI 53706-1388, USA}
\Email{\href{mailto:terwilli@math.wisc.edu}{terwilli@math.wisc.edu}}

\ArticleDates{Received December 22, 2012, in f\/inal form July 07, 2013; Published online July 15, 2013}

\Abstract{Let $\mathbb F$ denote a~f\/ield, and f\/ix a~nonzero $q\in\mathbb F$ such that $q^4\not=1$.
The universal Askey--Wilson algebra $\Delta_q$ is the associative $\mathbb F$-algebra def\/ined by
generators and relations in the following way.
The generators are $A$, $B$, $C$.
The relations assert that each of
\begin{gather*}
A+\frac{qBC-q^{-1}CB}{q^2-q^{-2}},
\qquad
B+\frac{qCA-q^{-1}AC}{q^2-q^{-2}},
\qquad
C+\frac{qAB-q^{-1}BA}{q^2-q^{-2}}
\end{gather*}
is central in $\Delta_q$.
The universal DAHA $\hat H_q$ of type $(C_1^\vee,C_1)$ is the associative $\mathbb F$-algebra def\/ined by
generators $\lbrace t^{\pm1}_i\rbrace_{i=0}^3$ and relations (i)~$t_i t^{-1}_i=t^{-1}_i t_i=1$; (ii)~$t_i+t^{-1}_i$ is central; (iii)~$t_0t_1t_2t_3=q^{-1}$.
We display an injection of $\mathbb F$-algebras $\psi:\Delta_q\to\hat H_q$ that sends
\begin{gather*}
A\mapsto t_1t_0+(t_1t_0)^{-1},
\qquad
B\mapsto t_3t_0+(t_3t_0)^{-1},
\qquad
C\mapsto t_2t_0+(t_2t_0)^{-1}.
\end{gather*}
For the map $\psi$ we compute the image of the three central elements mentioned above.
The algebra $\Delta_q$ has another central element of interest, called the Casimir element $\Omega$.
We compute the image of $\Omega$ under $\psi$.
We describe how the Artin braid group $B_3$ acts on $\Delta_q$ and $\hat H_q$ as a~group of automorphisms.
We show that $\psi$ commutes with these $B_3$ actions.
Some related results are obtained.}


\Classification{33D80; 33D45}

\section{Introduction}\label{Section1}

The Askey--Wilson polynomials were introduced in~\cite{awpoly} and soon became a~major topic in special
functions~\cite{ismail,koeswa}.
This topic became linked to representation theory through the work of A.~Zhedanov~\cite{Zhidd}.
In that article Zhedanov introduced the Askey--Wilson algebra ${\rm AW}(3)$, and showed that its ``ladder''
representations give the Askey--Wilson polynomials.
The algebra ${\rm AW}(3)$ is noncommutative and inf\/inite-dimensional.
It is def\/ined by generators and relations.
The relations involve a~scalar parameter $q$ and a~handful of extra scalar parameters.
The number of extra parameters ranges from~3 to~7 depending on which normalization is
used~\cite[equation~(6.1)]{koro}, \cite[Theorem~1.5]{aw}, \cite[Section~4.3]{wieg}, \cite[equations~(1.1a)--(1.1c)]{Zhidd}.
In~\cite{uaw} we introduced a~central extension of ${\rm AW}(3)$, denoted $\Delta_q$ and called the universal
Askey--Wilson algebra.
Up to normalization $\Delta_q$ is obtained from ${\rm AW}(3)$ by interpreting the extra parameters as central
elements in the algebra.
By construction $\Delta_q$ has no scalar parameters besides $q$, and there is a~surjective algebra
homomorphism $\Delta_q\to {\rm AW}(3)$.
One advantage of $\Delta_q$ over ${\rm AW}(3)$ is that $\Delta_q$ has a~larger automorphism group.
Our def\/inition of $\Delta_q$ was inspired by~\cite[Section~3]{dahater}, which in turn was
motivated by~\cite{ion}.

In~\cite{uaw} we began a~comprehensive investigation of $\Delta_q$.
In that paper we focused on its ring-theoretic aspects, and in a~followup paper~\cite{uawe} we related
$\Delta_q$ to the quantum algebra $U_q(\mathfrak{sl}_2)$.
In the present paper we relate $\Delta_q$ to the universal DAHA of type $(C_1^{\vee},C_1)$~\cite{dahater}.
This topic can be viewed as a~universal analog of a~topic considered by Koornwinder~\cite{Koo1, Koo2},
concerning how ${\rm AW}(3)$ is related to the DAHA of type $(C_1^{\vee},C_1)$.
We view~\cite{Koo1, Koo2} as groundbreaking and the main inspiration for the present paper.
In Section~\ref{Section16} we describe in detail how our results relate to those of Koornwinder~\cite{Koo1, Koo2}.

We will state our main results after we summarize the contents of~\cite{uaw,uawe}.

Our conventions for the paper are as follows.
An algebra is meant to be associative and have a~1.
A subalgebra has the same~1 as the parent algebra.
Fix a~f\/ield $\mathbb F$.
All unadorned tensor products are meant to be over $\mathbb F$.
We f\/ix a~nonzero $q\in\mathbb F$ such that $q^4\not=1$.
Recall the natural numbers $\mathbb N=\lbrace0,1,2,\ldots\rbrace$ and integers $\mathbb
Z=\lbrace0,\pm1,\pm2,\ldots\rbrace$.

The universal Askey--Wilson algebra $\Delta_q$ is the $\mathbb F$-algebra def\/ined by generators
and relations in the following way.
The generators are~$A$,~$B$,~$C$.
The relations assert that each of
\begin{gather}
A+\frac{qBC-q^{-1}CB}{q^2-q^{-2}},
\qquad
B+\frac{qCA-q^{-1}AC}{q^2-q^{-2}},
\qquad
C+\frac{qAB-q^{-1}BA}{q^2-q^{-2}}
\label{eq:introlist}
\end{gather}
is central in $\Delta_q$.
For the central elements~\eqref{eq:introlist} multiply each by $q+q^{-1}$ to get $\alpha$, $\beta$,
$\gamma$.
In~\cite{uaw} we obtained the following results about $\Delta_q$.
We gave an alternate presentation for $\Delta_q$ by generators and relations; the generators are~$A$,~$B$,~$\gamma$.
We gave a~faithful action of the modular group ${\rm{PSL}}_2(\mathbb Z)$ on $\Delta_q$ as a~group of
automorphisms; one generator sends $(A,B,C)\mapsto(B,C,A)$ and another generator sends
$(A,B,\gamma)\mapsto(B,A,\gamma)$.
We showed that $\lbrace A^iB^jC^k\alpha^r\beta^s\gamma^t\,|\,i,j,k,r,s,t\in\mathbb N\rbrace$ is a~basis for the
$\mathbb F$-vector space~$\Delta_q$.
We showed that the center $Z(\Delta_q)$ contains a~Casimir element
\begin{gather*}
\Omega=q^{-1}ACB+q^{-2}A^2+q^{-2}B^2+q^2C^2-q^{-1}A\alpha-q^{-1}B\beta-q C\gamma.
\end{gather*}
Under the assumption that $q$ is not a~root of unity, we showed that $Z(\Delta_q)$ is generated by
$\Omega$, $\alpha$, $\beta$, $\gamma$ and that $Z(\Delta_q)$ is isomorphic to a~polynomial algebra in four
variables.

In~\cite{uawe} we relate $\Delta_q$ to the quantum algebra $U_q(\mathfrak{sl}_2)$.
To describe this relationship we use the equitable presentation for $U_q(\mathfrak{sl}_2)$ which was
introduced in~\cite{equit}.
This equitable presentation has generators $x$, $y^{\pm1}$, $z$ and relations $y y^{-1}=y^{-1}y=1$,
\begin{gather*}
\frac{qxy-q^{-1}yx}{q-q^{-1}}=1,
\qquad
\frac{qyz-q^{-1}zy}{q-q^{-1}}=1,
\qquad
\frac{qzx-q^{-1}xz}{q-q^{-1}}=1.
\end{gather*}
Let $a$, $b$, $c$ denote mutually commuting indeterminates.
Let $\mathbb F\lbrack a^{\pm1},b^{\pm1},c^{\pm1}\rbrack$ denote the $\mathbb F$-algebra of Laurent
polynomials in $a$, $b$, $c$ that have all coef\/f\/icients in $\mathbb F$.
In~\cite[Theorems~2.16,~2.18]{uawe} we displayed an injection of $\mathbb F$-algebras $\natural:\Delta_q\to
U_q(\mathfrak{sl}_2)\otimes\mathbb F\lbrack a^{\pm1},b^{\pm1},c^{\pm1}\rbrack$ that sends
\begin{gather*}
A\mapsto x\otimes a+y\otimes a^{-1}+\frac{x y-y x}{q-q^{-1}}\otimes b c^{-1},
\qquad
B\mapsto y\otimes b+z\otimes b^{-1}+\frac{y z-z y}{q-q^{-1}}\otimes c a^{-1},
\\
C\mapsto z\otimes c+x\otimes c^{-1}+\frac{z x-x z}{q-q^{-1}}\otimes a b^{-1}.
\end{gather*}
The map $\natural$ sends
\begin{gather*}
\alpha\mapsto\Lambda\otimes\big(a+a^{-1}\big)+1\otimes\big(b+b^{-1}\big)\big(c+c^{-1}\big),
\\
\beta\mapsto\Lambda\otimes\big(b+b^{-1}\big)+1\otimes\big(c+c^{-1}\big)\big(a+a^{-1}\big),
\\
\gamma\mapsto\Lambda\otimes\big(c+c^{-1}\big)+1\otimes\big(a+a^{-1}\big)\big(b+b^{-1}\big),
\end{gather*}
where $\Lambda=qx+q^{-1}y+qz-qxyz$ is the normalized Casimir element of
$U_q(\mathfrak{sl}_2)$~\cite[Section~2]{uawe}.
In~\cite[Theorem~2.17]{uawe} we showed that $\natural$ sends $\Omega$ to{\samepage
\begin{gather*}
1\otimes\big(q+q^{-1}\big)^2-\Lambda^2\otimes1-1\otimes\big(a+a^{-1}\big)^2-1
\otimes\big(b+b^{-1}\big)^2-1\otimes\big(c+c^{-1}\big)^2
\\
\qquad
{}-\Lambda\otimes\big(a+a^{-1}\big)\big(b+b^{-1}\big)\big(c+c^{-1}\big).
\end{gather*}}

We now summarize the present paper.
We f\/irst show that the following is a~basis for the $\mathbb F$-vector space $\Delta_q$:
\begin{gather}
A^i C^j B^k\Omega^\ell\alpha^r\beta^s\gamma^t,
\qquad
j\in\lbrace0,1\rbrace,
\qquad
i,k,\ell,r,s,t\in\mathbb N.
\label{eq:introbase}
\end{gather}
This basis plays a~role in our main topic, which is about how $\Delta_q$ is related to the universal DAHA
$\hat H_q$ of type $(C_1^\vee,C_1)$.
The algebra $\hat H_q$ is a~variation on an algebra $\hat H$ introduced in~\cite{dahater}.
By def\/inition $\hat H_q$ is the $\mathbb F$-algebra with generators $\lbrace t^{\pm1}_i\rbrace_{i=0}^3$
and relations (i)~$t_i t^{-1}_i=t^{-1}_i t_i=1$; (ii)~$t_i+t^{-1}_i$ is central; (iii)~$t_0t_1t_2t_3=q^{-1}$.
We display an injection of $\mathbb F$-algebras $\psi:\Delta_q\to\hat H_q$ that sends
\begin{gather*}
A\mapsto t_1t_0+(t_1t_0)^{-1},
\qquad
B\mapsto t_3t_0+(t_3t_0)^{-1},
\qquad
C\mapsto t_2t_0+(t_2t_0)^{-1}.
\end{gather*}
The map $\psi$ sends
\begin{gather*}
\alpha\mapsto\big(q^{-1}t_0+q t^{-1}_0\big)\big(t_1+t^{-1}_1\big)+\big(t_2+t^{-1}_2\big)\big(t_3+t^{-1}_3\big),
\\
\beta\mapsto\big(q^{-1}t_0+q t^{-1}_0\big)\big(t_3+t^{-1}_3\big)+\big(t_1+t^{-1}_1\big)\big(t_2+t^{-1}_2\big),
\\
\gamma\mapsto\big(q^{-1}t_0+q t^{-1}_0\big)\big(t_2+t^{-1}_2\big)+\big(t_3+t^{-1}_3\big)\big(t_1+t^{-1}_1\big).
\end{gather*}
We show that $\psi$ sends $\Omega$ to
\begin{gather*}
\big(q+q^{-1}\big)^2-(q^{-1}t_0+qt_0^{-1})^2-\big(t_1+t^{-1}_1\big)^2-\big(t_2+t^{-1}_2\big)^2-\big(t_3+t^{-1}_3\big)^2
\\
\qquad
{}-\big(q^{-1}t_0+q t_0^{-1}\big)\big(t_1+t^{-1}_1\big)\big(t_2+t^{-1}_2\big)\big(t_3+t^{-1}_3\big).
\end{gather*}
We remark that for the above results some parts are easier to prove than others.
It is relatively easy to show that~$\psi$ exists as an algebra homomorphism.
Indeed this existence essentially follows from~\cite[Theorem~5.2]{dahater}, although in our formal argument
we take another approach which quickly yields the result from f\/irst principles.
We found it relatively hard to show that $\psi$ is injective; indeed this argument takes up the majority of
the paper.
To establish injectivity we display a~basis for~$\hat H_q$, and use it to show that~$\psi$ sends the
basis~\eqref{eq:introbase} to a~linearly independent set.
Adapting~\cite[Theorem~2.6]{ion},~\cite[Lemma~4.2]{dahater} we show how the Artin braid group~$B_3$ acts on~$\hat H_q$ as a~group of automorphisms.
The group $B_3$ is a~homomorphic preimage of~${\rm{PSL}}_2(\mathbb Z)$, and we mentioned earlier that
${\rm{PSL}}_2(\mathbb Z)$ acts on~$\Delta_q$ as a~group of automorphisms.
Pulling back this ${\rm{PSL}}_2(\mathbb Z)$ action we get a~$B_3$ action on~$\Delta_q$ as a~group of
automorphisms.
We show that~$\psi$ commutes with the~$B_3$ actions for~$\Delta_q$ and~$\hat H_q$.
Now consider the image of~$\Delta_q$ under~$\psi$.
Adapting~\cite[Theorem~5.1]{Koo2} we show that the subalgebra $\lbrace h\in\hat H_q|t_0h=h t_0\rbrace$ is
generated by this image together with~$t_0$ and $\lbrace t_i+t^{-1}_i\rbrace_{i=1}^3$.
For this subalgebra we give a~presentation by generators and relations.
Roughly speaking, this presentation amounts to a~$q$-analog of~\cite[Theorem~2.1]{oblomkov} and a~universal
analog of~\cite[Def\/inition~6.1, Corollary~6.3]{Koo1}.
Under the assumption that $q$ is not a~root of unity, we show that~$Z(\hat H_q)$ is generated by $\lbrace
t_i+t^{-1}_i\rbrace_{i=0}^3$ and that~$Z(\hat H_q)$ is isomorphic to a~polynomial algebra in four variables.
Roughly speaking, this is a~universal analog of~\cite[Theorem~5.3]{Koo2}.

The paper is organized as follows.
In Section~\ref{Section2}, after reviewing $\Delta_q$ we obtain a~basis for this algebra that will be useful.
In Section~\ref{Section3} we def\/ine $\hat H_q$ and discuss its symmetries.
In Section~\ref{Section4} we state f\/ive theorems which describe an injection $\psi:\Delta_q\to\hat H_q$; these
are Theorems~\ref{thm:main1}--\ref{thm:main3}.
In Section~\ref{Section5} we establish some identities in $\hat H_q$ that will be used repeatedly.
In Section~\ref{Section6} we prove Theorems~\ref{thm:main1},~\ref{thm:main4},~\ref{thm:main4a}.
In Section~\ref{Section7} we display a~basis for $\hat H_q$, along with some reduction rules that show how to
write any given element of $\hat H_q$ in the basis.
Sections~\ref{Section8},~\ref{Section9} are devoted to proving Theorem~\ref{thm:main2}.
Sections~\ref{Section10}--\ref{Section12} are devoted to proving Theorem~\ref{thm:main3}.
In Section~\ref{Section13} we consider the image of $\Delta_q$ under the map $\psi$.
We show that the subalgebra $\lbrace h\in\hat H_q|t_0h=h t_0\rbrace$ is generated by this image together
with $t_0$ and $\lbrace t_i+t^{-1}_i\rbrace_{i=1}^3$.
In Section~\ref{Section14} we give a~presentation for this subalgebra by generators and relations.
In Section~\ref{Section15} we describe $Z(\hat H_q)$.
In Section~\ref{Section16} we compare our results with those of Koornwinder~\cite{Koo1,Koo2}.

\section[The universal Askey-Wilson algebra]{The universal Askey--Wilson algebra}\label{Section2}

We now begin our formal argument.
In this section we discuss the universal Askey--Wilson algebra.
After reviewing its basic features we establish a~basis for the algebra that will be useful later in the
paper.
\begin{Definition}[\protect{\cite[Def\/inition~1.2]{uaw}}] \label{def:uaw}
Def\/ine an $\mathbb F$-algebra $\Delta_q$ by generators and relations in the following way.
The generators are~$A$,~$B$,~$C$.
The relations assert that each of
\begin{gather}
A+\frac{qBC-q^{-1}CB}{q^2-q^{-2}},
\qquad
B+\frac{qCA-q^{-1}AC}{q^2-q^{-2}},
\qquad
C+\frac{qAB-q^{-1}BA}{q^2-q^{-2}}
\label{eq:comlist}
\end{gather}
is central in $\Delta_q$.
The algebra $\Delta_q$ is called the {\it universal Askey--Wilson algebra}.
\end{Definition}

\begin{Definition}[\protect{\cite[Def\/inition~1.3]{uaw}}]
\label{def:abc}
For the three central elements in~\eqref{eq:comlist}, multiply each by $q+q^{-1}$ to get
$\alpha$, $\beta$, $\gamma$.
Thus
\begin{gather}
A+\frac{qBC-q^{-1}CB}{q^2-q^{-2}}=\frac{\alpha}{q+q^{-1}},
\label{eq:u1}
\\
B+\frac{qCA-q^{-1}AC}{q^2-q^{-2}}=\frac{\beta}{q+q^{-1}},
\label{eq:u2}
\\
C+\frac{qAB-q^{-1}BA}{q^2-q^{-2}}=\frac{\gamma}{q+q^{-1}}.
\label{eq:u3}
\end{gather}
Note that each of $\alpha$, $\beta$, $\gamma$ is central in $\Delta_q$.
Note also that $A$, $B$, $\gamma$ is a~generating set for $\Delta_q$.
\end{Definition}

We now discuss some automorphisms of $\Delta_q$.
Recall that the modular group ${\rm{PSL}}_2(\mathbb Z)$ has a~presentation by generators~$p$,~$s$ and
relations $p^3=1$, $s^2=1$.
See for example~\cite{alpern}.
\begin{Lemma}[\protect{\cite[Theorem~3.1]{uaw}}]
 The group ${\rm{PSL}}_2(\mathbb Z)$ acts on $\Delta_q$ as a~group of
automorphisms in the following way:

\centerline{
\begin{tabular}
{c| ccc | c c c}
{}$u$ & $A$ & $B$ & $C$ & $\alpha$ & $\beta$ & $\gamma$
\\[0.5mm]
\hline
{}$p(u)$ & $B$ & $C$ & $A$ & $\beta$ & $\gamma$ & $\alpha$
\\[0.5mm]
{}$s(u)$ & $B$ & $A$ & $C+\dfrac{AB-BA}{q-q^{-1}}$ & $\beta$ & $\alpha$ & $\gamma$
\end{tabular}
}
\end{Lemma}

The group ${\rm{PSL}}_2(\mathbb Z)$ has a~central extension called the Artin braid group $B_3$.
The group $B_3$ is def\/ined as follows.

\begin{Definition}
\label{def:braidgroup}
The Artin braid group $B_3$ is def\/ined by generators $\rho$, $\sigma$ and relations $\rho^3=\sigma^2$.
For notational convenience def\/ine $\tau=\rho^3=\sigma^2$.
\end{Definition}

There exists a~group homomorphism $B_3\to{\rm{PSL}}_2(\mathbb Z)$ that sends $\rho\mapsto p$ and
$\sigma\mapsto s$.
Via this homomorphism we pull back the ${\rm{PSL}}_2(\mathbb Z)$ action on $\Delta_q$, to get a~$B_3$
action on $\Delta_q$ as a~group of automorphisms.
This action is described as follows.
\begin{Lemma}\label{thm:try}
The group $B_3$ acts on $\Delta_q$ as a~group of automorphisms such that $\tau(h)=h$ for all $h\in\Delta_q$
and $\rho$, $\sigma$ do the following:

\centerline{
\begin{tabular}
{c| ccc | c c c}
{}$u$ & $A$ & $B$ & $C$ & $\alpha$ & $\beta$ & $\gamma$
\\[0.5mm]
\hline
{}$\rho(u)$ & $B$ & $C$ & $A$ & $\beta$ & $\gamma$ & $\alpha$
\\[0.5mm]
{}$\sigma(u)$ & $B$ & $A$ & $C+\dfrac{AB-BA}{q-q^{-1}}$ & $\beta$ & $\alpha$ & $\gamma$
\end{tabular}
}
\end{Lemma}

In Def\/inition~\ref{def:abc} we def\/ined the central elements $\alpha$, $\beta$, $\gamma$ of~$\Delta_q$.
There is another central element of interest, called the Casimir element~$\Omega$.
This element is def\/ined as follows.
\begin{Definition}[\protect{\cite[Lemma~6.1]{uaw}}]
\label{def:casdelta}
Def\/ine an element $\Omega\in\Delta_q$ by
\begin{gather}
\Omega=q^{-1}ACB+q^{-2}A^2+q^{-2}B^2+q^{2}C^2-q^{-1}A\alpha-q^{-1}B\beta-q C\gamma.
\label{eq:casdelta}
\end{gather}
We call $\Omega$ the {\it Casimir element} of $\Delta_q$.
\end{Definition}

\begin{Lemma}[\protect{\cite[Theorems~6.2,~8.2]{uaw}}] The Casimir element $\Omega$ is contained in the center $Z(\Delta_q)$.
Moreover $\lbrace\Omega^i\alpha^r\beta^s\gamma^t\, |\, i,r,s,t\in\mathbb N\rbrace$ is a~basis for the $\mathbb
F$-vector space $Z(\Delta_q)$, provided that $q$ is not a~root of unity.
\end{Lemma}

\begin{Lemma}[\protect{\cite[Theorem~6.4]{uaw}}] 
 The Casimir element $\Omega$ is f\/ixed by everything in $B_3$.
\end{Lemma}

Given an $\mathbb F$-algebra $\mathcal A$, by an {\it antiautomorphism} of $\mathcal A$ we mean
an $\mathbb F$-linear bijection $\zeta:\mathcal A\to\mathcal A$ such that $(uv)^\zeta=v^\zeta u^\zeta$ for
all $u,v\in\mathcal A$.
\begin{Lemma}
\label{lem:daggeranti}
There exists an antiautomorphism $\dagger$ of $\Delta_q$ that sends
\begin{gather*}
A\mapsto B,
\qquad
B\mapsto A,
\qquad
C\mapsto C,
\qquad
\alpha\mapsto\beta,
\qquad
\beta\mapsto\alpha,
\qquad
\gamma\mapsto\gamma.
\end{gather*}
Moreover $\dagger^2=1$.
\end{Lemma}

\begin{proof}
Routine using~\eqref{eq:u1}--\eqref{eq:u3}.
\end{proof}

\begin{Lemma}
The Casimir element $\Omega$ is fixed by the antiautomorphism $\dagger$ from Lemma~{\rm \ref{lem:daggeranti}}.
\end{Lemma}
\begin{proof}
This follows from~\cite[Lemma~6.1]{uaw}.
\end{proof}

We mention how $\Delta_q$ and $\Delta_{q^{-1}}$ are related.
\begin{Lemma}
\label{lem:isomap}
There exists an isomorphism of $\mathbb F$-algebras $\xi:\Delta_q\to\Delta_{q^{-1}}$ that sends
\begin{gather*}
A\mapsto B,
\qquad
B\mapsto A,
\qquad
C\mapsto C,
\qquad
\alpha\mapsto\beta,
\qquad
\beta\mapsto\alpha,
\qquad
\gamma\mapsto\gamma.
\end{gather*}
\end{Lemma}
\begin{proof}
Routine using~\eqref{eq:u1}--\eqref{eq:u3}.
\end{proof}

\begin{Lemma}
The isomorphism $\xi:\Delta_q\to\Delta_{q^{-1}}$ from Lemma~{\rm \ref{lem:isomap}} sends the Casimir element of
$\Delta_q$ to the Casimir element of $\Delta_{q^{-1}}$.
\end{Lemma}

\begin{proof}
This follows from~\cite[Lemma~6.1]{uaw}.
\end{proof}

We are going to display a~basis for the $\mathbb F$-vector space $\Delta_q$.
Two such bases can be found in~\cite[Theorems~4.1,~7.5]{uaw}, but these are not suited for our present
purpose.
To obtain a~suitable basis we work with the following presentation of $\Delta_q$.
\begin{Proposition}
\label{prop:thirdpres}
The $\mathbb F$-algebra $\Delta_q$ is presented by generators and relations in the following way.
The generators are $A$, $B$, $C$, $\Omega$, $\alpha$, $\beta$, $\gamma$.
The relations assert that each of $\Omega$, $\alpha$, $\beta$, $\gamma$ is central and
\begin{gather*}
BA = q^2AB+q\big(q^2-q^{-2}\big)C-q\big(q-q^{-1}\big)\gamma,
\\
BC = q^{-2}CB-q^{-1}\big(q^2-q^{-2}\big)A+q^{-1}\big(q-q^{-1}\big)\alpha,
\\
CA = q^{-2}AC-q^{-1}(q^{2}-q^{-2})B+q^{-1}\big(q-q^{-1}\big)\beta,
\\
C^2 = q^{-2}\Omega-q^{-3}ACB-q^{-4}A^2-q^{-4}B^2+q^{-3}A\alpha+q^{-3}B\beta+q^{-1}C\gamma.
\end{gather*}
\end{Proposition}

\begin{proof}
Referring to the above four equations, the f\/irst three are reformulations of~\eqref{eq:u1}--\eqref{eq:u3}
and the fourth is a~reformulation of~\eqref{eq:casdelta}.
\end{proof}
\begin{Definition}
The generators $A$, $B$, $C$, $\Omega$, $\alpha$, $\beta$, $\gamma$ of $\Delta_q$ are called {\it balanced}.
\end{Definition}
\begin{Note}
\label{note:third}
Referring to the presentation of $\Delta_q$ from Proposition~\ref{prop:thirdpres}, consider the
relations which assert that $\Omega$, $\alpha$, $\beta$, $\gamma$ are central.
These relations can be expressed as
\begin{alignat*}{4}
&\Omega A=A\Omega,
\qquad &&
\Omega B=B\Omega,
\qquad&&
\Omega C=C\Omega,&
\\
&\alpha A=A\alpha,
\qquad&&
\alpha B=B\alpha,
\qquad&&
\alpha C=C\alpha,&
\\
&\beta A=A\beta,
\qquad&&
\beta B=B\beta,
\qquad&&
\beta C=C\beta,&
\\
&\gamma A=A\gamma,
\qquad&&
\gamma B=B\gamma,
\qquad&&
\gamma C=C\gamma,&
\\
&\alpha\Omega=\Omega\alpha,
\qquad&&
\beta\Omega=\Omega\beta,
\qquad&&
\gamma\Omega=\Omega\gamma,&
\\
&\beta\alpha=\alpha\beta,
\qquad&&
\gamma\alpha=\alpha\gamma,
\qquad&&
\gamma\beta=\beta\gamma.&
\end{alignat*}
\end{Note}
\begin{Definition}
By a~{\it reduction rule} for $\Delta_q$ we mean an equation which appears in
Proposition~\ref{prop:thirdpres} or Note~\ref{note:third}.
A~reduction rule from Proposition~\ref{prop:thirdpres} is said to be of the {\it first kind}.
A~reduction rule from Note~\ref{note:third} is said to be of the {\it second kind}.
\end{Definition}

\begin{Definition}
 For an integer $n\geq0$, by a~{\it word of length $n$} in $\Delta_q$ we mean a~product $g_1g_2\cdots
g_n$ such that $g_i$ is a~balanced generator of $\Delta_q$ for $1\leq i\leq n$.
We interpret the word of length 0 as the multiplicative identity in $\Delta_q$.
A word is called {\it forbidden} whenever it is the left-hand side of a~reduction rule.
Every forbidden word has length two.
A forbidden word is said to be of the {\it first kind} (resp.\
{\it second kind}) whenever the corresponding reduction rule is of the f\/irst (resp.\
second) kind.
\end{Definition}

\begin{Definition}
 Let $w$ denote a~forbidden word in $\Delta_q$, and consider the corresponding reduction rule.
By a~{\it descendent of~$w$} we mean a~word that appears on the right-hand side of that reduction rule.
\end{Definition}

\begin{Example}
The descendents of $BA$ are $AB$, $C$, $\gamma$.
The descendents of $BC$ are $CB$, $A$, $\alpha$.
The descendents of $CA$ are $AC$, $B$, $\beta$.
The descendents of $C^2$ are $\Omega$, $ACB$, $A^2$, $B^2$, $A\alpha$, $B\beta$,~$C\gamma$.
A forbidden word of the second kind has a~single descendent, obtained by interchanging its two factors.
\end{Example}

\begin{Theorem}
\label{thm:usefulbasis}
The following is a~basis for the $\mathbb F$-vector space $\Delta_q$:
\begin{gather}
A^i C^j B^k\Omega^\ell\alpha^r\beta^s\gamma^t,
\qquad
j\in\lbrace0,1\rbrace,
\qquad
i,k,\ell,r,s,t\in\mathbb N.
\label{eq:usefulbasis}
\end{gather}
\end{Theorem}

\begin{proof}
We invoke Bergman's Diamond Lemma~\cite[Theorem~1.2]{berg}.
Let $g_1g_2\cdots g_n$ denote a~word in $\Delta_q$.
This word is called {\it reducible} whenever there exists an integer $i$ $(2\leq i\leq n)$ such that
$g_{i-1}g_i$ is forbidden.
The word is called {\it irreducible} whenever it is not reducible.
The list~\eqref{eq:usefulbasis} consists of the irreducible words in $\Delta_q$.
Let $w=g_1g_2\cdots g_n$ denote a~word in $\Delta_q$.
By an {\it inversion} in $w$ we mean an ordered pair of integers $(i,j)$ such that $1\leq i<j\leq n$ and
the word $g_i g_j$ is forbidden.
The inversion $(i,j)$ is of the {\it first kind} (resp.\
{\it second kind}) whenever the forbidden word $g_ig_j$ is of the f\/irst kind (resp.\
second kind).
Let $W$ denote the set of all words in $\Delta_q$.
We def\/ine a~partial order $<$ on $W$ as follows.
Pick any words $w$, $w'$ in $W$ and write $w=g_1g_2\cdots g_n$.
We say that $w$ {\it dominates} $w'$ whenever there exists an integer $i$ $(2\leq i\leq n)$ such that
$(i-1,i)$ is an inversion for $w$, and $w'$ is obtained from $w$ by replacing $g_{i-1}g_i$ by one of its
descendents.
In this case either (i)~$w$~has more inversions of the f\/irst kind than $w'$, or (ii)~$w$~and~$w'$ have
the same number of inversions of the f\/irst kind, but $w$ has more inversions of the second kind than~$w'$.
By these comments the transitive closure of the domination relation on $W$ is a~partial order on $W$ which
we denote by $<$.
By construction $<$ is a~semigroup partial order~\cite[p.~181]{berg} and satisf\/ies the descending chain
condition~\cite[p.~179]{berg}.
We now relate the partial order $<$ to our reduction rules.
Let $w=g_1g_2\cdots g_n$ denote a~reducible word in $\Delta_q$.
Then there exists an integer $i$ $(2\leq i\leq n)$ such that $g_{i-1}g_i$ is forbidden.
There exists a~reduction rule with $g_{i-1}g_i$ on the left-hand side; in $w$ we eliminate $g_{i-1}g_i$
using this reduction rule and thereby express $w$ as a~linear combination of words, each less than $w$ with
respect to~$<$.
Therefore the reduction rules are compatible with $<$ in the sense of Bergman~\cite[p.~181]{berg}.
In order to employ the Diamond Lemma, we must show that the ambiguities are resolvable in the sense of
Bergman~\cite[p.~181]{berg}.
There are potentially two kinds of ambiguities; inclusion ambiguities and overlap
ambiguities~\cite[p.~181]{berg}.
For the present example there are no inclusion ambiguities.
The nontrivial overlap ambiguities are
\begin{gather*}
BCA,
\qquad
BC^2,
\qquad
C^2A.
\end{gather*}
Take for example $BCA$.
The words $BC$ and $CA$ are forbidden.
Therefore $BCA$ can be reduced in two ways; we could evaluate $BC$ f\/irst or we could evaluate~$CA$
f\/irst.
Either way, after a~4-step reduction we get the same resolution, which is
\begin{gather*}
q^{-3}\big(q^2-q^{-2}\big)\Omega+q^{-6}ACB-q^{-3}\big(q^4-q^{-4}\big)A^2-q^{-3}\big(q^4-q^{-4}\big)B^2
\\
\qquad
{} + q^{-3}\big(q^3-q^{-3}\big)A\alpha+q^{-3}\big(q^3-q^{-3}\big)B\beta+q^{-3}\big(q-q^{-1}\big)C\gamma.
\end{gather*}
Therefore the ambiguity $BCA$ is resolvable.
The ambiguities $BC^2$ and $C^2A$ are similarly shown to be resolvable.
The resolution of $BC^2$ is
\begin{gather*}
q^{-6}B\Omega-q^{-7}ACB^2-q^{-8}A^2B-q^{-8}B^3+q^{-7}AB\alpha+q^{-7}B^2\beta+q^{-5}CB\gamma
\\
\qquad
{}- q^{-3}\big(q^4-q^{-4}\big)AC+q^{-2}\big(q^2-q^{-2}\big)C\alpha
+q^{-4}\big(q^2-q^{-2}\big)^2B\\
\qquad
{}
-q^{-4}\big(q-q^{-1}\big)\big(q^2-q^{-2}\big)\beta
\end{gather*}
and the resolution of $C^2A$ is
\begin{gather*}
q^{-6}A\Omega-q^{-7}A^2CB-q^{-8}AB^2-q^{-8}A^3+q^{-7}A^2\alpha+q^{-7}AB\beta+q^{-5}AC\gamma
\\
\qquad
{}- q^{-3}\big(q^4-q^{-4}\big)CB+q^{-2}\big(q^2-q^{-2}\big)C\beta
+q^{-4}\big(q^2-q^{-2}\big)^2A\\
\qquad
{}
-q^{-4}\big(q-q^{-1}\big)\big(q^2-q^{-2}\big)\alpha.
\end{gather*}
We conclude that every ambiguity is resolvable, so by the Diamond Lemma~\cite[Theorem~1.2]{berg}
the irreducible words form a~basis for~$\Delta_q$.
The result follows.
\end{proof}

\section[The universal DAHA $\hat H_q$ of type $(C_1^\vee,C_1)$]{The universal
DAHA $\boldsymbol{\hat H_q}$ of type $\boldsymbol{(C_1^\vee,C_1)}$}\label{Section3}

The double af\/f\/ine Hecke algebra (DAHA) for a~reduced root system was def\/ined by
Cherednik~\cite{chered}, and the def\/inition was extended to include nonreduced root systems by
Sahi~\cite{sahi}.
The most general DAHA of rank 1 is associated with the root system $(C_1^\vee,C_1)$.
In~\cite{dahater} we introduced a~central extension of this algebra called the universal DAHA of type
$(C_1^\vee,C_1)$.
In the present paper we will work with a~variation on this algebra.

For notational convenience def\/ine a~four element set
\begin{gather*}
\mathbb I=\lbrace0,1,2,3\rbrace.
\end{gather*}

The following def\/inition is a~variation on~\cite[Def\/inition~3.1]{dahater}.
\begin{Definition}
\label{def:udaha1}
Let $\hat H_q$ denote the $\mathbb F$-algebra def\/ined by generators $\lbrace
t^{\pm1}_i\rbrace_{i\in\mathbb I}$ and relations
\begin{gather}
t_i t^{-1}_i=t^{-1}_i t_i=1, \qquad
i\in\mathbb I,
\nonumber
\\
t_i+t^{-1}_i
\ \text{is central},
\qquad
i\in\mathbb I,
\label{eq:centralsum}
\\
t_0t_1t_2t_3=q^{-1}.
\nonumber
\end{gather}
We call $\hat H_q$ the {\it universal DAHA of type $(C^\vee_1,C_1)$}.
\end{Definition}

\begin{Remark}
In~\cite[Def\/inition~3.1]{dahater} we def\/ined an $\mathbb F$-algebra $\hat H$ by generators $\lbrace
t^{\pm1}_i\rbrace_{i\in\mathbb I}$ and relations (i)~$t_i t^{-1}_i=t^{-1}_i t_i=1$; (ii)~$t_i+t^{-1}_i$ is
central; (iii)~$t_0t_1t_2t_3$ is central.
The algebra $\hat H_q$ is a~homomorphic image of $\hat H$.
\end{Remark}

The following two lemmas are immediate from Def\/inition~\ref{def:udaha1}.
\begin{Lemma}
In the algebra $\hat H_q$ the scalar $q^{-1}$ is equal to each of the following:
\begin{gather*}
t_0t_1t_2t_3,
\qquad
t_1t_2t_3t_0,
\qquad
t_2t_3t_0t_1,
\qquad
t_3t_0t_1t_2.
\end{gather*}
\end{Lemma}
\begin{Lemma}
\label{lem:z4aut}
There exists an automorphism of $\hat H_q$ that sends $t_0\mapsto t_1\mapsto t_2\mapsto t_3\mapsto t_0$.
\end{Lemma}

Recall the braid group $B_3$ from Def\/inition~\ref{def:braidgroup}.
\begin{Lemma}
\label{thm:b31}
The group $B_3$ acts on $\hat H_q$ as a~group of automorphisms such that $\tau(h)=t^{-1}_0ht_0$ for all
$h\in{\hat H_q}$ and $\rho$, $\sigma$ do the following:

 \centerline{
\begin{tabular}{c|cccc}
{}$h$ & $t_0$ & $t_1$ & $t_2$ & $t_3$\bsep{1pt}
\\[0.5mm]
\hline
{}$\rho(h)$ & $t_0$ & $t^{-1}_0t_3t_0$ & $t_1$ & $t_2$\tsep{2pt}
\\[0.5mm]
{}$\sigma(h)$ & $t_0$ & $t^{-1}_0t_3t_0$ & $t_1t_2t^{-1}_1$ & $t_1$
\end{tabular}
}
\end{Lemma}
\begin{proof}
This is routinely checked, or see~\cite[Lemma~4.2]{dahater}.
\end{proof}
\begin{Lemma}\label{thm:half}
The $B_3$ action on $\hat H_q$ does the following to the central elements~{\rm \eqref{eq:centralsum}}.
The generator $\tau$ fixes every central element.
The ge\-ne\-rators $\rho$, $\sigma$ satisfy the table below.

 \centerline{
\begin{tabular}{c|cccc}
{}$h$ & $t_0+t^{-1}_0$ & $t_1+t^{-1}_1$ & $t_2+t^{-1}_2$ & $t_3+t^{-1}_3$\bsep{1pt}
\\
\hline
{}$\rho(h)$ & $t_0+t^{-1}_0$ & $t_3+t^{-1}_3$ & $t_1+t^{-1}_1$ & $t_2+t^{-1}_2$\tsep{2pt}
\\[0.5mm]
{}$\sigma(h)$ & $t_0+t^{-1}_0$ & $t_3+t^{-1}_3$ & $t_2+t^{-1}_2$ & $t_1+t^{-1}_1$
\end{tabular}
}
\end{Lemma}
\begin{proof}
Use~\eqref{eq:centralsum} and Lemma~\ref{thm:b31}.
\end{proof}

\begin{Lemma}
\label{lem:Hantiaut}
There exists a~unique antiautomorphism $\dagger$ of $\hat H_q$ that sends
\begin{gather*}
t_0\mapsto t_0,
\qquad
t_1\mapsto t_3,
\qquad
t_2\mapsto t_2,
\qquad
t_3\mapsto t_1.
\end{gather*}
Moreover $\dagger^2=1$.
\end{Lemma}

\begin{proof}
Use Def\/inition~\ref{def:udaha1}.
\end{proof}

\begin{Lemma}
\label{lem:Hqiso}
There exists a~unique isomorphism of $\mathbb F$-algebras $\xi:\hat H_q\to\hat H_{q^{-1}}$ that sends
\begin{gather*}
t_0\mapsto t^{-1}_0,
\qquad
t_1\mapsto t^{-1}_3,
\qquad
t_2\mapsto t^{-1}_2,
\qquad
t_3\mapsto t^{-1}_1.
\end{gather*}
\end{Lemma}

\begin{proof}
Use Def\/inition~\ref{def:udaha1}.
\end{proof}

\section[How $\Delta_q$ is related to $\hat H_q$]{How $\boldsymbol{\Delta_q}$
is related to $\boldsymbol{\hat H_q}$}\label{Section4}

In this section we state f\/ive theorems concerning how $\Delta_q$ is related to $\hat H_q$.
The proofs of these theorems will take up most of the rest of the paper.
\begin{Theorem}
\label{thm:main1}
There exists a~unique $\mathbb F$-algebra homomorphism $\psi:\Delta_q\to\hat H_q$ that sends
\begin{gather*}
A\mapsto t_1t_0+(t_1t_0)^{-1},
\qquad
B\mapsto t_3t_0+(t_3t_0)^{-1},
\qquad
C\mapsto t_2t_0+(t_2t_0)^{-1}.
\end{gather*}

The homomorphism $\psi$ sends
\begin{gather*}
\alpha\mapsto\big(q^{-1}t_0+q t^{-1}_0\big)\big(t_1+t^{-1}_1\big)+\big(t_2+t^{-1}_2\big)\big(t_3+t^{-1}_3\big),
\\
\beta\mapsto\big(q^{-1}t_0+q t^{-1}_0\big)\big(t_3+t^{-1}_3\big)+\big(t_1+t^{-1}_1\big)\big(t_2+t^{-1}_2\big),
\\
\gamma\mapsto\big(q^{-1}t_0+q t^{-1}_0\big)\big(t_2+t^{-1}_2\big)+\big(t_3+t^{-1}_3\big)\big(t_1+t^{-1}_1\big).
\end{gather*}
\end{Theorem}

\begin{Theorem}
\label{thm:main4}
For all $g\in B_3$ the following diagram commutes:
\begin{gather*}
\begin{CD}
\Delta_q@>\psi>>\hat H_q
\\
@V g VV@VV g V
\\
\Delta_q@>>\psi>\hat H_q
\end{CD}
\end{gather*}
\end{Theorem}

\begin{Theorem}
\label{thm:main4a}
The following diagrams commute:
\begin{gather*}
\begin{CD}
\Delta_q@>\psi>>\hat H_q
\\
@V\dagger VV@VV\dagger V
\\
\Delta_q@>>\psi>\hat H_q
\end{CD}
\qquad
\begin{CD}
\Delta_q@>\psi>>\hat H_q
\\
@V\xi VV@VV\xi V
\\
\Delta_{q^{-1}}@>>\psi>\hat H_{q^{-1}}
\end{CD}
\end{gather*}
\end{Theorem}
\begin{Theorem}
\label{thm:main2}
Under the homomorphism $\psi$ from Theorem~{\rm \ref{thm:main1}} the image of $\Omega$ is
\begin{gather}
\big(q+q^{-1}\big)^2-\big(q^{-1}t_0+qt_0^{-1}\big)^2-\big(t_1+t^{-1}_1\big)^2-\big(t_2+t^{-1}_2\big)^2
-\big(t_3+t^{-1}_3\big)^2
\nonumber
\\
\qquad
{}-\big(q^{-1}t_0+q t_0^{-1}\big)\big(t_1+t^{-1}_1\big)\big(t_2+t^{-1}_2\big)\big(t_3+t^{-1}_3\big).
\label{eq:omimage}
\end{gather}
\end{Theorem}

\begin{Theorem}\label{thm:main3}
The homomorphism $\psi$ from Theorem~{\rm \ref{thm:main1}} is injective.
\end{Theorem}

\section[Preliminaries concerning $\hat H_q$]{Preliminaries concerning $\boldsymbol{\hat H_q}$}\label{Section5}

In this section we establish some basic facts about~$\hat H_q$.
These facts will be used repeatedly for the rest of the paper.
\begin{Definition}
For the algebra $\hat H_q$ def\/ine
\begin{gather}
T_i=t_i+t^{-1}_i,
\qquad
i\in\mathbb I.
\label{eq:Ti}
\end{gather}
Note that each $T_i$ is central in $\hat H_q$.
\end{Definition}

In Def\/inition~\ref{def:udaha1} we gave a~presentation for $\hat H_q$ involving the generators
$\lbrace t^{\pm1}_i\rbrace_{i\in\mathbb I}$.
Sometimes it is convenient to use $\lbrace T_i\rbrace_{i\in\mathbb I}$ instead of $\lbrace
t^{-1}_i\rbrace_{i\in\mathbb I}$.
In terms of the generators $\lbrace t_i\rbrace_{i\in\mathbb I}$, $\lbrace T_i\rbrace_{i\in\mathbb I}$ the
algebra $\hat H_q$ looks as follows.
\begin{Lemma}
\label{prop:altpres}
The $\mathbb F$-algebra $\hat H_q$ has a~presentation by generators $\lbrace t_i\rbrace_{i\in\mathbb
I}$, $\lbrace T_i\rbrace_{i\in\mathbb I}$ and relations
\begin{gather*}
t^2_i=T_i t_i-1,
\qquad
i\in\mathbb I,
\\
T_i
\
{\mbox{is central}},
\qquad
i\in\mathbb I,
\\
t_0t_1t_2t_3=q^{-1}.
\end{gather*}
\end{Lemma}

\begin{Definition}
 Let $X$, $Y$ denote the following elements of $\hat H_q$:
\begin{gather}
X=t_3t_0,
\qquad
Y=t_0t_1.
\label{eq:XY}
\end{gather}
Note that each of $X$, $Y$ is invertible.
\end{Definition}

\begin{Lemma}
\label{lem:123}
For the algebra $\hat H_q$,
\begin{gather}
\label{eq:t123}
t_1=t^{-1}_0Y,
\qquad
t_2=q^{-1}Y^{-1}t_0X^{-1},
\qquad
t_3=Xt^{-1}_0.
\end{gather}
Moreover $\hat H_q$ is generated by $X^{\pm1}$, $Y^{\pm1}$, $t^{\pm1}_0$.
\end{Lemma}

\begin{proof}
The relations~\eqref{eq:t123} are routinely checked using Def\/inition~\ref{def:udaha1} and~\eqref{eq:XY}.
\end{proof}

In terms of the generators $X^{\pm1}$, $Y^{\pm1}$, $t^{\pm1}_0$ the $\lbrace T_i\rbrace_{i\in\mathbb
I}$ look as follows.
\begin{Lemma}
\label{lem:T123}
For the algebra $\hat H_q$ the following $(i)$--$(iv)$ hold.
\begin{enumerate}\itemsep=0pt
\item[$(i)$] $T_0=t_0+t^{-1}_0$.
\item[$(ii)$] $T_1$ is equal to each of
\begin{gather*}
t^{-1}_0Y+Y^{-1}t_0,
\qquad
Yt^{-1}_0+t_0Y^{-1}.
\end{gather*}
\item[$(iii)$] $T_2$ is equal to each of
\begin{gather*}
qt^{-1}_0YX+q^{-1}X^{-1}Y^{-1}t_0,
\qquad
q Xt^{-1}_0Y+q^{-1}Y^{-1}t_0X^{-1},
\\
q YXt^{-1}_0+q^{-1}t_0X^{-1}Y^{-1}.
\end{gather*}
\item[$(iv)$] $T_3$ is equal to each of
\begin{gather*}
t^{-1}_0X+X^{-1}t_0,
\qquad
Xt^{-1}_0+t_0X^{-1}.
\end{gather*}
\end{enumerate}
\end{Lemma}

\begin{proof}
$(i)$ Clear.

$(ii)$ Using the equation on the left in~\eqref{eq:t123},
\begin{gather*}
T_1=t_1+t^{-1}_1=t^{-1}_0Y+Y^{-1}t_0.
\end{gather*}
Also
\begin{gather*}
T_1=YT_1Y^{-1}=Y\big(t^{-1}_0Y+Y^{-1}t_0\big)Y^{-1}=Yt^{-1}_0+t_0Y^{-1}.
\end{gather*}

$(iii)$ Using the middle equation in~\eqref{eq:t123},
\begin{gather*}
T_2=t^{-1}_2+t_2=q Xt^{-1}_0Y+q^{-1}Y^{-1}t_0X^{-1}.
\end{gather*}
Also
\begin{gather*}
T_2=X^{-1}T_2X=X^{-1}\big(q Xt^{-1}_0Y+q^{-1}Y^{-1}t_0X^{-1}\big)X=q t^{-1}_0YX+q^{-1}X^{-1}Y^{-1}t_0
\end{gather*}
and
\begin{gather*}
T_2=YT_2Y^{-1}=Y\big(q Xt^{-1}_0Y+q^{-1}Y^{-1}t_0X^{-1}\big)Y^{-1}=q YXt^{-1}_0+q^{-1}t_0X^{-1}Y^{-1}.
\end{gather*}

$(iv)$ Similar to the proof of $(ii)$ above.
\end{proof}

In Section~\ref{Section3} we discussed some automorphisms and antiautomorphisms of $\hat H_q$.
We now consider how these maps act on~$X$,~$Y$.
The following four lemmas are routinely checked.

\begin{Lemma}
\label{lem:z4}
Consider the automorphism of $\hat H_q$ from Lemma~{\rm \ref{lem:z4aut}}.
This automorphism sends
\begin{gather*}
X\mapsto Y\mapsto q^{-1}X^{-1}\mapsto q^{-1}Y^{-1}\mapsto X.
\end{gather*}
\end{Lemma}
\begin{Lemma}
Consider the automorphisms $\rho$, $\sigma$ of $\hat H_q$ from Lemma~{\rm \ref{thm:b31}}.
The automorhism $\rho$ sends
\begin{gather*}
X\mapsto q^{-1}Y^{-1}t_0X^{-1}t_0,
\qquad
Y\mapsto X.
\end{gather*}
The automorphism $\sigma$ sends
\begin{gather*}
X\mapsto t_0^{-1}Y t_0,
\qquad
Y\mapsto X.
\end{gather*}
\end{Lemma}

\begin{Lemma}
Recall the antiautomorphism $\dagger$ of $\hat H_q$ from Lemma~{\rm \ref{lem:Hantiaut}}.
This map swaps $X$, $Y$.
\end{Lemma}
\begin{Lemma}
Recall the isomorphism $\xi:\hat H_q\to\hat H_{q^{-1}}$ from Lemma~{\rm \ref{lem:Hqiso}}.
This map sends $X\mapsto Y^{-1}$ and $Y\mapsto X^{-1}$.
\end{Lemma}

We now give some relations that show how $t_0$ commutes past the $X^{\pm1}$, $Y^{\pm1}$.
\begin{Lemma}
\label{lem:4rel}
The following relations hold in $\hat H_q$:
\begin{gather}
t_0X=X^{-1}t_0+X T_0-T_3,
\label{eq:4rel1}
\\
t_0X^{-1}=Xt_0-XT_0+T_3,
\label{eq:4rel2}
\\
t_0Y=Y^{-1}t_0+Y T_0-T_1,
\label{eq:4rel3}
\\
t_0Y^{-1}=Yt_0-Y T_0+T_1.
\label{eq:4rel4}
\end{gather}
\end{Lemma}
\begin{proof}
To obtain~\eqref{eq:4rel1},~\eqref{eq:4rel2} replace $t_0^{-1}$ by $T_0-t_0$ in Lemma~\ref{lem:T123}$(iv)$.
To obtain~\eqref{eq:4rel3},~\eqref{eq:4rel4} replace $t_0^{-1}$ by $T_0-t_0$ in Lemma~\ref{lem:T123}$(ii)$.
\end{proof}

We now consider how $X$, $Y$ are related.
\begin{Lemma}
\label{lem:tprod1}
The following relations hold in $\hat H_q$:
\begin{alignat}{3}\label{eq:line2}
&t_0t_2=q^{-1}t^{-1}_3T_1-q^{-1}Y X^{-1},
\qquad&&
t_0^{-1}t^{-1}_2=q t_1T_3-q X^{-1}Y,&
\\
\label{eq:line2t}
&t_1t_3=q^{-1}t^{-1}_0T_2-q^{-2}X^{-1}Y^{-1},
\qquad&&
t_1^{-1}t^{-1}_3=q t_2T_0-Y^{-1}X^{-1},&
\\
\label{eq:line2tt}
&t_2t_0=q^{-1}t^{-1}_1T_3-q^{-1}Y^{-1}X,
\qquad&&
t_2^{-1}t^{-1}_0=q t_3T_1-qXY^{-1},&
\\
\label{eq:line2ttt}
&t_3t_1=q^{-1}t^{-1}_2T_0-X Y,
\qquad&&
 t_3^{-1}t^{-1}_1=q t_0T_2-q^2YX.&
\end{alignat}
\end{Lemma}

\begin{proof}
Concerning~\eqref{eq:line2}, the equation on the left comes from
\begin{gather*}
q^{-1}YX^{-1}=t_0t_1^2t_2=t_0(T_1t_1-1)t_2=q^{-1}t_3^{-1}T_1-t_0t_2.
\end{gather*}
The equation on the right comes from
\begin{gather*}
q X^{-1}Y=t^{-1}_0t_3^{-2}t^{-1}_2=t^{-1}_0\big(T_3t^{-1}_3-1\big)t^{-1}_2=q t_1T_3-t^{-1}_0t^{-1}_2.
\end{gather*}
To obtain~\eqref{eq:line2t}--\eqref{eq:line2ttt}, repeatedly apply the automorphism from
Lemma~\ref{lem:z4aut} to everything in~\eqref{eq:line2}, and use Lemma~\ref{lem:z4}.
\end{proof}
\begin{Definition}
\label{def:ci}
Let $\lbrace C_i\rbrace_{i\in\mathbb I}$ denote the following elements in $\hat H_q$:
\begin{gather*}
C_0=q\big(qYX-q^{-1}XY\big),
\\
C_1=-\big(q^{-1}YX^{-1}-qX^{-1}Y\big),
\\
C_2=q^{-1}\big(qY^{-1}X^{-1}-q^{-1}X^{-1}Y^{-1}\big),
\\
C_3=-\big(q^{-1}Y^{-1}X-qXY^{-1}\big).
\end{gather*}
\end{Definition}
\begin{Lemma}
\label{lem:cadjust}
The automorphism from Lemma~{\rm \ref{lem:z4aut}} sends
\begin{gather*}
C_0\mapsto C_1\mapsto C_2\mapsto C_3\mapsto C_0.
\end{gather*}
\end{Lemma}
\begin{proof}
Use Lemma~\ref{lem:z4} and Def\/inition~\ref{def:ci}.
\end{proof}

\begin{Proposition}
\label{prop:Ciform}
The following relations hold in $\hat H_q$:
\begin{gather}
\label{eq:C0}
C_0=qT_2t_0+T_3t_1+q^{-1}T_0t_2+T_1t_3-q^{-1}T_0T_2-T_1T_3,
\\
\label{eq:C1}
C_1=T_2t_0+qT_3t_1+T_0t_2+q^{-1}T_1t_3-T_0T_2-q^{-1}T_1T_3,
\\
\label{eq:C2}
C_2=q^{-1}T_2t_0+T_3t_1+qT_0t_2+T_1t_3-q^{-1}T_0T_2-T_1T_3,
\\
\label{eq:C3}
C_3=T_2t_0+q^{-1}T_3t_1+T_0t_2+q T_1t_3-T_0T_2-q^{-1}T_1T_3.
\end{gather}
\end{Proposition}
\begin{proof}
To verify~\eqref{eq:C0}, use~\eqref{eq:line2ttt} together with
\begin{gather*}
t_3^{-1}t_1^{-1}=(T_3-t_3)(T_1-t_1)
=T_1T_3-t_3T_1-t_1T_3+t_3t_1.
\end{gather*}
To verify~\eqref{eq:C1}--\eqref{eq:C3}, repeatedly apply the automorphism from Lemma~\ref{lem:z4aut} to
everything in~\eqref{eq:C0}, and use Lemma~\ref{lem:cadjust}.
\end{proof}

We mention a~result for future use.
\begin{Lemma}
\label{lem:sigma02}
The automorphism $\sigma$ of $\hat H_q$ sends
\begin{alignat*}{3}
&t_1t_3\mapsto q^{-1}t_0^{-1}t_2^{-1},
\qquad&&
t_3^{-1}t_1^{-1}\mapsto q t_2t_0,&
\\
&t_0t_2\mapsto q^{-1}t_3^{-1}t_1^{-1},
\qquad&&
t_2^{-1}t_0^{-1}\mapsto q t_1t_3.&
\end{alignat*}
\end{Lemma}
\begin{proof}
This is routinely checked using the action of $\sigma$ given in Lemma~\ref{thm:b31}.
\end{proof}

\section{The proof of Theorems~\ref{thm:main1},~\ref{thm:main4},~\ref{thm:main4a}}\label{Section6}

In this section we prove the f\/irst three theorems from Section~\ref{Section4}.
\begin{Lemma}[\protect{\cite[Lemma~3.8]{dahater}}]
\label{lem:why1}
  Let $u$, $v$ denote invertible elements in any algebra such that each of
$u+u^{-1}$, $v+v^{-1}$ is central.
Then
\begin{enumerate}\itemsep=0pt
\item[$(i)$] $uv+(uv)^{-1}=vu+(vu)^{-1}$; \item[$(ii)$] $uv+(uv)^{-1}$ commutes with each of $u$, $v$.
\end{enumerate}
\end{Lemma}

\begin{proof}
$(i)$ Write $U=u+u^{-1}$ and $V=v+v^{-1}$.
We have both
\begin{gather*}
uv+(uv)^{-1}=uv+(V-v)(U-u)=uv+vu-vU-uV+UV,
\\
vu+(vu)^{-1}=vu+(U-u)(V-v)=vu+uv-uV-vU+UV.
\end{gather*}
The result follows.

$(ii)$ Using $(i)$ we have
\begin{gather*}
u^{-1}\big(uv+v^{-1}u^{-1}\big)u=vu+u^{-1}v^{-1}=uv+v^{-1}u^{-1}.
\end{gather*}
Therefore $uv+(uv)^{-1}$ commutes with $u$.
Similarly $uv+(uv)^{-1}$ commutes with $v$.
\end{proof}
\begin{Corollary}
\label{cor:titj}
For distinct $i,j\in\mathbb I$,
\begin{enumerate}\itemsep=0pt
\item[$(i)$] $t_it_j+(t_it_j)^{-1}=t_jt_i+(t_jt_i)^{-1}$.
\item[$(ii)$] $t_it_j+(t_it_j)^{-1}$ commutes with each of $t_i$, $t_j$.
\end{enumerate}
\end{Corollary}

\begin{proof}
By Lemma~\ref{lem:why1} and since $t_k+t^{-1}_k$ is central for $k\in\mathbb I$.
\end{proof}
\begin{Definition}
\label{def:xyz}
We def\/ine elements $A$, $B$, $C$ in $\hat H_q$ as follows:
\begin{gather}
A=t_1t_0+(t_1t_0)^{-1}=t_0t_1+(t_0t_1)^{-1}=Y+Y^{-1},
\nonumber
\\
B=t_3t_0+(t_3t_0)^{-1}=t_0t_3+(t_0t_3)^{-1}=X+X^{-1},
\nonumber
\\
C=t_2t_0+(t_2t_0)^{-1}=t_0t_2+(t_0t_2)^{-1}.
\label{eq:xyz3}
\end{gather}
\end{Definition}

\begin{Lemma}
\label{cor:txyz}
In the algebra $\hat H_q$ the element $t_0$ commutes with each of $A$, $B$, $C$.
\end{Lemma}
\begin{proof}
By Corollary~\ref{cor:titj}$(ii)$ and Def\/inition~\ref{def:xyz}.
\end{proof}

The following is a~variation on~\cite[Theorem~5.1]{dahater}.
\begin{Lemma}
\label{lem:B3xyz}
The $B_3$ action on $\hat H_q$ does the following to the elements $A$, $B$, $C$ from Definition~{\rm \ref{def:xyz}}.
The generator $\tau$ fixes each of $A$, $B$, $C$.
The generator $\rho$ sends $A\mapsto B\mapsto C\mapsto A$.
The generator $\sigma$ swaps $A$, $B$ and sends $C\mapsto C'$ where
\begin{gather*}
q C+q^{-1}C'+AB=q^{-1}C+q C'+BA
=\big(q^{-1}t_0+q t^{-1}_0\big)T_2+T_1T_3.
\end{gather*}
\end{Lemma}

\begin{proof}
The generator $\tau$ f\/ixes each of $A$, $B$, $C$ by Lemma~\ref{cor:txyz} and since $\tau(h)=t^{-1}_0h t_0$
for all $h\in\hat H_q$.
The generator $\rho$ sends $A\mapsto B\mapsto C\mapsto A$ by Lemma~\ref{thm:b31} and
Def\/inition~\ref{def:xyz}.
Similarly the generator $\sigma$ swaps $A$, $B$.
Def\/ine $C'=\sigma(C)$.
We show that $C'$ satisf\/ies the equations of the lemma statement.
We f\/irst show that
\begin{gather}
\label{eq:B3part1}
qC+q^{-1}C'+AB=\big(q^{-1}t_0+q t^{-1}_0\big)T_2+T_1T_3.
\end{gather}
Since $A=Y+Y^{-1}$ and $B=X+X^{-1}$,
\begin{gather}
AB=YX+YX^{-1}+Y^{-1}X+Y^{-1}X^{-1}.
\label{eq:xyexpand}
\end{gather}
By~\eqref{eq:xyz3} along with~\eqref{eq:line2} and~\eqref{eq:line2tt},
\begin{gather*}
C=t_0t_2+t^{-1}_2t_0^{-1}
=\big(qt_3+q^{-1}t^{-1}_3\big)T_1-qXY^{-1}-q^{-1}YX^{-1}.
\end{gather*}
Using~\eqref{eq:xyz3} and Lemma~\ref{lem:sigma02} along with~\eqref{eq:line2t} and~\eqref{eq:line2ttt},
\begin{gather}
C'=qt_1t_3+q^{-1}t^{-1}_3t^{-1}_1
=T_0T_2-qYX-q^{-1}X^{-1}Y^{-1}.
\label{eq:zpexpand}
\end{gather}
To verify~\eqref{eq:B3part1}, evaluate the left-hand side using~\eqref{eq:xyexpand}--\eqref{eq:zpexpand}
and simplify the result using Def\/inition~\ref{def:ci}, Proposition~\ref{prop:Ciform}, and~\eqref{eq:Ti}.
We have verif\/ied~\eqref{eq:B3part1}.
Next we show that
\begin{gather}
\label{eq:B3part2}
q^{-1}C+q C'+BA=\big(q^{-1}t_0+q t^{-1}_0\big)T_2+T_1T_3.
\end{gather}
To obtain~\eqref{eq:B3part2}, apply $\sigma$ to each side of~\eqref{eq:B3part1} and evaluate the result.
To aid in this evaluation, recall that $\sigma$ swaps $A$, $B$; also $\sigma$ swaps $C$, $C'$ since
$\sigma^2=\tau$ and $\tau(C)=C$.
By these comments and Lemma~\ref{thm:half} we routinely obtain~\eqref{eq:B3part2}.
\end{proof}

The following is a~variation on~\cite[Theorem~5.2]{dahater}.
\begin{Proposition}
\label{thm:Qlevel2}
In the algebra $\hat H_q$ the elements $A$, $B$, $C$ are related as follows:
\begin{gather*}
A+\frac{qBC-q^{-1}CB}{q^2-q^{-2}}=\frac{\big(q^{-1}t_0+qt^{-1}_0\big)T_1+T_2T_3}{q+q^{-1}},
\\
B+\frac{qCA-q^{-1}AC}{q^2-q^{-2}}=\frac{\big(q^{-1}t_0+qt^{-1}_0\big)T_3+T_1T_2}{q+q^{-1}},
\\
C+\frac{qAB-q^{-1}BA}{q^2-q^{-2}}=\frac{\big(q^{-1}t_0+qt^{-1}_0\big)T_2+T_3T_1}{q+q^{-1}}.
\end{gather*}
\end{Proposition}
\begin{proof}
To get the last equation, eliminate $C'$ from the equations of Lemma~\ref{lem:B3xyz}.
To get the other two equations use the $B_3$ action from Lemma~\ref{thm:b31}.
Specif\/ically, apply $\rho$ twice to the last equation and use the data in Lemma~\ref{thm:half}, together
with the fact that $\rho$ cyclically permutes~$A$, $B$, $C$ and f\/ixes~$t_0$.
\end{proof}

\begin{proof}[Proof of Theorem~\ref{thm:main1}] Immediate from Lemma~\ref{cor:txyz} and
Proposition~\ref{thm:Qlevel2}.
\end{proof}

Back in Def\/inition~\ref{def:abc} we def\/ined some elements $\alpha$, $\beta$, $\gamma$ of
$\Delta_q$.
From now on we retain the notation $\alpha$, $\beta$, $\gamma$ for their images under the map
$\psi:\Delta_q\to\hat H_q$.
Thus the elements $\alpha$, $\beta$, $\gamma$ of $\hat H_q$ satisfy
\begin{gather}
\alpha=\big(q^{-1}t_0+qt_0^{-1}\big)T_1+T_2T_3,
\label{eq:alphainH}
\\
\beta=\big(q^{-1}t_0+qt_0^{-1}\big)T_3+T_1T_2,
 \label{eq:betainH}
\\
\gamma=\big(q^{-1}t_0+qt_0^{-1}\big)T_2+T_3T_1.
\label{eq:gammainH}
\end{gather}

\begin{proof}[Proof of Theorem~\ref{thm:main4}] Without loss we may assume $g=\rho$ or $g=\sigma$.
By Lemma~\ref{thm:try} the action of~$\rho$ on $\Delta_q$ cyclically permutes $A$, $B$, $C$.
By Lemma~\ref{lem:B3xyz} the action of~$\rho$ on $\hat H_q$ cyclically permutes $A$, $B$, $C$.
By Lemma~\ref{thm:try} the action of~$\sigma$ on~$\Delta_q$ swaps $A$, $B$ and f\/ixes $\gamma$.
The action of~$\sigma$ on~$\hat H_q$ swaps $A$, $B$ by Lemma~\ref{lem:B3xyz}.
The action of~$\sigma$ on~$\hat H_q$ f\/ixes $\gamma$ by~\eqref{eq:gammainH} and
Lemmas~\ref{thm:b31},~\ref{thm:half}.
The result follows.
\end{proof}

\begin{proof}[Proof of Theorem~\ref{thm:main4a}] In each case, chase $A$, $B$, $C$ around the diagram using
Theorem~\ref{thm:main1} and Corollary~\ref{cor:titj}$(i)$, together with Lemma~\ref{lem:daggeranti}
and~\ref{lem:Hantiaut} for $\dagger$ and with Lemma~\ref{lem:isomap} and~\ref{lem:Hqiso} for~$\xi$.
\end{proof}

\section[A basis for the $\mathbb F$-vector space $\hat H_q$]{A basis
for the $\boldsymbol{\mathbb F}$-vector space $\boldsymbol{\hat H_q}$}\label{Section7}

Our next general goal is to prove Theorem~\ref{thm:main2}.
The proof will be completed in Section~\ref{Section9}.
In the present section we obtain a~basis for the $\mathbb F$-vector space $\hat H_q$.
The basis consists of
\begin{gather}
Y^i X^j t_0^k T_0^\ell T_1^r T_2^s T_3^t,
\qquad
i,j\in\mathbb Z,
\qquad
k\in\lbrace0,1\rbrace,
\qquad
\ell,r,s,t\in\mathbb N.
\label{eq:basisxy1pre}
\end{gather}
We also obtain a~set of relations for $\hat H_q$ called reduction rules.
The reduction rules show how to write any given element of $\hat H_q$ as a~linear combination of the basis
elements~\eqref{eq:basisxy1pre}.

To begin the basis project, we are going to display a~presentation of $\hat H_q$ that
contains detailed information about how the generators commute past each other.
We will give two versions of this presentation.
For version I we attempt to optimize clarity.
For version II we attempt to optimize utility.
We hope that taken together the two versions are reasonably clear and useful.
The relations in version II become our reduction rules.

We now give version I.
\begin{Proposition}
\label{prop:nicelong}
The $\mathbb F$-algebra $\hat H_q$ has a~presentation by generators $X^{\pm1}$, $Y^{\pm1}$, $\lbrace
t_i\rbrace_{i\in\mathbb I}$, $\lbrace T_i\rbrace_{i\in\mathbb I}$, $\lbrace C_i\rbrace_{i\in\mathbb I}$ and
relations $XX^{-1}=1$, $X^{-1}X=1$, $YY^{-1}=1$, $Y^{-1}Y=1$, the $\lbrace T_i\rbrace_{i\in\mathbb I}$ are
central,
\begin{gather*}
t_0^2=t_0T_0-1,
\\
t_1=(T_0-t_0)Y,
\\
t_2=q^{-1}Y^{-1}t_0X^{-1},
\\
t_3=X(T_0-t_0),
\\
t_0X=X^{-1}t_0+X T_0-T_3,
\\
t_0X^{-1}=Xt_0-XT_0+T_3,
\\
t_0Y=Y^{-1}t_0+Y T_0-T_1,
\\
t_0Y^{-1}=Y t_0-Y T_0+T_1,
\\
XY=q^2YX-C_0,
\\
X^{-1}Y=q^{-2}YX^{-1}+q^{-1}C_1,
\\
X^{-1}Y^{-1}=q^2Y^{-1}X^{-1}-q^2C_2,
\\
XY^{-1}=q^{-2}Y^{-1}X+q^{-1}C_3,
\\
C_0=qT_2t_0+T_3t_1+q^{-1}T_0t_2+T_1t_3-q^{-1}T_0T_2-T_1T_3,
\\
C_1=T_2t_0+qT_3t_1+T_0t_2+q^{-1}T_1t_3-T_0T_2-q^{-1}T_1T_3,
\\
C_2=q^{-1}T_2t_0+T_3t_1+qT_0t_2+T_1t_3-q^{-1}T_0T_2-T_1T_3,
\\
C_3=T_2t_0+q^{-1}T_3t_1+T_0t_2+q T_1t_3-T_0T_2-q^{-1}T_1T_3.
\end{gather*}
\end{Proposition}

\begin{proof}
Consider the relations in the proposition statement.
We now show that these relations hold in $\hat H_q$.
This is clear for the relations shown in the line, so consider the 16 displayed relations.
Displayed relation~1 is from Lemma~\ref{prop:altpres}.
Displayed relations~2--4 follow from Lemma~\ref{lem:123}.
Displayed relations 5--8 are from Lemma~\ref{lem:4rel}.
Displayed relations 9--12 are from Def\/inition~\ref{def:ci}.
Displayed relations 13--16 are from Proposition~\ref{prop:Ciform}.
We have shown that the relations in the proposition statement hold in $\hat H_q$.
Conversely, one routinely checks that the relations in the proposition statement imply the def\/ining
relations for $\hat H_q$ given in Lemma~\ref{prop:altpres}.
\end{proof}

We now give version II.
Roughly speaking, this version amounts to a~universal analog of~\cite[Proposition~5.2]{Koo1}.
\begin{Proposition}
\label{prop:hgenrel}
The $\mathbb F$-algebra $\hat H_q$ has a~presentation by generators $X^{\pm1}$, $Y^{\pm1}$, $t_0$, $\lbrace
T_i\rbrace_{i\in\mathbb I}$ and relations $XX^{-1}=1$, $X^{-1}X=1$, $YY^{-1}=1$, $Y^{-1}Y=1$, the $\lbrace
T_i\rbrace_{i\in\mathbb I}$ are central,
\begin{gather*}
t_0^2=t_0T_0-1,
\\
t_0X=X^{-1}t_0+X T_0-T_3,
\\
t_0X^{-1}=X t_0-X T_0+T_3,
\\
t_0Y=Y^{-1}t_0+Y T_0-T_1,
\\
t_0Y^{-1}=Y t_0-Y T_0+T_1,
\\
XY=q^2YX-q t_0T_2+q^{-1}T_0T_2+Y^{-1}t_0T_3-q^{-2}Y^{-1}T_0T_3
\\
\hphantom{XY=}{}
 + q^{-2}Y^{-1}X T_0^2-q^{-2}Y^{-1}X t_0T_0-X T_0T_1+X t_0T_1,
\\
X^{-1}Y=q^{-2}YX^{-1}+\big(q-q^{-1}\big)q^{-1}T_1T_3-q^{-1}T_0T_2+q^{-1}t_0T_2-Y^{-1}t_0T_3
\\
\hphantom{X^{-1}Y=}{}
+ q^{-2}X T_0T_1-q^{-2}X t_0T_1+q^{-2}Y^{-1}T_0T_3-q^{-2}Y^{-1}X T_0^2+q^{-2}Y^{-1}X t_0T_0,
\\
X^{-1}Y^{-1}=q^2Y^{-1}X^{-1}-q^2Y^{-1}T_0T_3+q^2Y^{-1}t_0T_3+q T_0T_2-q t_0T_2
\\
\hphantom{X^{-1}Y^{-1}=}{}
- q^2X T_0T_1+q^2X t_0T_1+q^2Y^{-1}X T_0^2-q^2Y^{-1}X t_0T_0,
\\
XY^{-1}=q^{-2}Y^{-1}X+X T_0T_1-X t_0T_1-q^{-2}Y^{-1}X T_0^2+q^{-2}Y^{-1}X t_0T_0
\\
\hphantom{XY^{-1}=}{}
+ q^{-2}Y^{-1}T_0T_3-q^{-2}Y^{-1}t_0T_3-q^{-1}T_0T_2+q^{-1}t_0T_2.
\end{gather*}
\end{Proposition}

\begin{proof}\sloppy
In Proposition~\ref{prop:nicelong} eliminate $\lbrace t_i\rbrace_{i=1}^3$ using the displayed relations
2--4, and eliminate $\lbrace C_i\rbrace_{i\in\mathbb I}$ using the displayed relations 13--16.
Simplify the results using the displayed relations \mbox{5--8}.
\end{proof}

We just gave two versions of a~presentation for $\hat H_q$.
From now on we focus on version II.
This version will yield our reduction rules and basis for $\hat H_q$.
\begin{Definition}
 The generators $X^{\pm1}$, $Y^{\pm1}$, $t_0$, $\lbrace T_i\rbrace_{i\in\mathbb I}$ of $\hat H_q$ are called {\it
balanced}.
\end{Definition}

\begin{Note}
\label{note:centralmean}
 Referring to the presentation of $\hat H_q$ from Proposition~\ref{prop:hgenrel}, consider the relations
which assert that the $\lbrace T_i\rbrace_{i\in\mathbb I}$ are central.
These relations can be expressed as
\begin{gather*}
T_i X^{\pm1}=X^{\pm1}T_i,
\qquad
T_i Y^{\pm1}=Y^{\pm1}T_i,
\qquad
T_i t_0=t_0T_i,
\qquad
i\in\mathbb I,
\\
T_i T_j=T_jT_i,
\qquad
i,j\in\mathbb I,
\qquad
i>j.
\end{gather*}
\end{Note}

\begin{Definition}
 By a~{\it reduction rule} for $\hat H_q$ we mean an equation that appears in
Proposition~\ref{prop:hgenrel} or Note~\ref{note:centralmean}.
Of these reduction rules, the last four in Proposition~\ref{prop:hgenrel} are said to be of the {\it first
kind}, the preceeding f\/ive are said to be of the {\it second kind}, and the rest are said to be of the
{\it third kind}.
\end{Definition}
\begin{Definition}
For an integer $n\geq0$, by a~{\it word of length $n$} in $\hat H_q$ we mean a~product $g_1g_2\cdots
g_n$ such that $g_i$ is a~balanced generator of $\hat H_q$ for $1\leq i\leq n$.
We interpret the word of length 0 as the multiplicative identity in $\hat H_q$.
A word is called {\it forbidden} whenever it is the left-hand side of a~reduction rule.
Every forbidden word has length two.
A forbidden word is said to be of the {\it first kind} (resp.\
{\it second kind}) (resp.\
{\it third kind}) whenever the corresponding reduction rule is of the f\/irst kind (resp.\
second kind) (resp.\
third kind).
\end{Definition}
\begin{Definition}
Let $w$ denote a~forbidden word in $\hat H_q$, and consider the corresponding reduction rule.
By a~{\it descendent of~$w$} we mean a~word that appears on the right-hand side of that reduction rule.
\end{Definition}

Roughly speaking, the following result amounts to a~universal analog of~\cite[Theorem~5.3]{Koo1}.
\begin{Proposition}
\label{prop:basisv1}
The following is a~basis for the $\mathbb F$-vector space $\hat H_q$:
\begin{gather}
Y^i X^j t_0^k T_0^\ell T_1^r T_2^s T_3^t
\qquad
i,j\in\mathbb Z,
\qquad
k\in\lbrace0,1\rbrace,
\qquad
\ell,r,s,t\in\mathbb N.
\label{eq:basisxy1}
\end{gather}
\end{Proposition}

\begin{proof}
We invoke Bergman's Diamond Lemma~\cite[Theorem~1.2]{berg}.
Let $g_1g_2\cdots g_n$ denote a~word in $\hat H_q$.
This word is called {\it reducible} whenever there exists an integer $i$ $(2\leq i\leq n)$ such that
$g_{i-1}g_i$ is forbidden.
A word is called {\it irreducible} whenever it is not reducible.
The list~\eqref{eq:basisxy1} consists of the irreducible words in $\hat H_q$.
Let $w=g_1g_2\dots g_n$ denote a~word in $\hat H_q$.
By an {\it inversion} in $w$ we mean an ordered pair of integers $(i,j)$ such that $1\leq i<j\leq n$ and
the word $g_i g_j$ is forbidden.
The inversion $(i,j)$ is of the {\it first kind} (resp.\
{\it second kind}) (resp.\
{\it third kind}) whenever the forbidden word~$g_ig_j$ is of the f\/irst kind (resp.\
second kind) (resp.\
third kind).
Let~$W$ denote the set of all words in $\hat H_q$.
We def\/ine a~partial order $<$ on $W$ as follows.
Pick any words $w$, $w'$ in $W$ and write $w=g_1g_2\cdots g_n$.
We say that~$w$ {\it dominates} $w'$ whenever there exists an integer~$i$ $(2\leq i\leq n)$ such that
$(i-1,i)$ is an inversion for $w$, and $w'$ is obtained from~$w$ by replacing $g_{i-1}g_i$ by one of its
descendents.
In this case either (i)~$w$~has more inversions of the f\/irst kind than~$w'$, or (ii)~$w$~and~$w'$ have
the same number of inversions of the f\/irst kind, but $w$ has more inversions of the second kind than~$w'$, or (iii)~$w$~and~$w'$ have the same number of inversions for each of the f\/irst and second kind, but~$w$ has more inversions of the third kind than~$w'$.
By these comments the transitive closure of the domination relation on~$W$ is a~partial order on~$W$ which
we denote by~$<$.
By construction $<$ is a~semigroup partial order~\cite[p.~181]{berg} and satisf\/ies the descending chain
condition~\cite[p.~179]{berg}.
We now relate the partial order $<$ to our reduction rules.
Let $w=g_1g_2\cdots g_n$ denote a~reducible word in~$\hat H_q$.
Then there exists an integer~$i$ $(2\leq i\leq n)$ such that $g_{i-1}g_i$ is forbidden.
There exists a~reduction rule with $g_{i-1}g_i$ on the left-hand side; in~$w$ we eliminate $g_{i-1}g_i$
using this reduction rule and thereby express $w$ as a~linear combination of words, each less than $w$ with
respect to~$<$.
Therefore the reduction rules are compatible with~$<$ in the sense of Bergman~\cite[p.~181]{berg}.
In order to employ the Diamond Lemma, we must show that the ambiguities are resolvable in the sense of
Bergman~\cite[p.~181]{berg}.
There are potentially two kinds of ambiguities; inclusion ambiguities and overlap
ambiguities~\cite[p.~181]{berg}.
For the present example there are no inclusion ambiguities.
The nontrivial overlap ambiguities are
\begin{gather*}
t_0X Y,
\qquad
t_0X^{-1}Y,
\qquad
t_0X Y^{-1},
\qquad
t_0X^{-1}Y^{-1},
\\
t_0^2X,
\qquad
t_0^2X^{-1},
\qquad
t_0^2Y,
\qquad
t_0^2Y^{-1},
\\
X X^{-1}Y,
\qquad
X X^{-1}Y^{-1},
\qquad
X^{-1}X Y,
\qquad
X^{-1}X Y^{-1},
\\
X YY^{-1},
\qquad
X Y^{-1}Y,
\qquad
X^{-1}Y Y^{-1},
\qquad
X^{-1}Y^{-1}Y,
\\
t_0XX^{-1},
\qquad
t_0X^{-1}X,
\qquad
t_0YY^{-1},
\qquad
t_0Y^{-1}Y.
\end{gather*}
Take for example $t_0X Y$.
The words $t_0X$ and $XY$ are forbidden.
Therefore $t_0XY$ can be reduced in two ways; we could evaluate $t_0X$ f\/irst or we could evaluate $XY$
f\/irst.
Either way, after a~3-step reduction we get the same resolution, which is
\begin{gather*}
q^2Y X T_0+q^{-2}YX^{-1}T_0+q^2Y^{-1}X^{-1}t_0+q^2Y^{-1}X T_0-q^2X T_1
\\
\qquad{}
+ \big(q^{-2}-1\big)XT_0^2T_1+\big(1-q^{-2}\big)X t_0T_0T_1-X^{-1}T_1-Y T_3-q^2Y^{-1}T_3
\\
\qquad{}
- \big(q-q^{-1}\big)t_0T_0T_2+\big(1-q^{-2}\big)T_0T_1T_3+q T_2.
\end{gather*}
Therefore the ambiguity $t_0X Y$ is resolvable.
The other ambiguities listed above are similarly shown to be resolvable.
Their resolutions are displayed in the tables below.

\vspace{1mm}

\centerline{
\begin{tabular}{c| c }
Ambiguity & Resolution
\\
\hline
{}$t_0X^{-1}Y$ & $q^{-2}Y^{-1}X t_0-q^{-2}Y^{-1}XT_0+Y T_3+q^{-2}Y^{-1}T_3-q^{-1}T_2$,\tsep{2pt}
\\[1.75mm]
{}$t_0X Y^{-1}$ & $q^{-2}Y X^{-1}t_0-q^{-2}Y X^{-1}T_0+(q^{-2}-1)X t_0T_0T_1+(1-q^{-2})X T_0^2T_1$
\\
& $
\qquad
+ \; q^{-2}X T_1 + X^{-1} T_1 + (1-q^{-2}) t_0 T_1 T_3 + (q^{-2}-1)T_0 T_1 T_3 - q^{-1} T_2$,
\\[1.75mm]
{}$t_0 X^{-1} Y^{-1}$&$q^2 Y X t_0 -q^2 Y X T_0 + q T_2$
\\[1.75mm]
{}$t_0^2X$ & $X^{-1}t_0T_0+X T_0^2-X-T_0T_3$
\\[1.75mm]
{}$t_0^2X^{-1}$ & $X t_0T_0-X T_0^2-X^{-1}+T_0T_3$
\\[1.75mm]
{}$t_0^2Y$ & $Y^{-1}t_0T_0+Y T_0^2-Y-T_0T_1$
\\[1.75mm]
{}$t_0^2Y^{-1}$ & $Y t_0T_0-Y T_0^2-Y^{-1}+T_0T_1$
\end{tabular}
}

\vspace{2mm}

\centerline{
\begin{tabular}{c| ccc c }
Ambiguity & $XX^{-1}Y$ & $XX^{-1}Y^{-1}$ & $X^{-1}XY$ & $X^{-1}XY^{-1}$
\\
\hline
Resolution & $Y$ & $Y^{-1}$ & $Y$ & $Y^{-1}$
\end{tabular}
}

\vspace{2mm}

\centerline{
\begin{tabular}{c| ccc c }
Ambiguity & $XYY^{-1}$ & $XY^{-1}Y$ & $X^{-1}YY^{-1}$ & $X^{-1}Y^{-1}Y$
\\
\hline
Resolution & $X$ & $X$ & $X^{-1}$ & $X^{-1}$
\end{tabular}
}

\vspace{2mm}

\centerline{
\begin{tabular}{c| ccc c }
Ambiguity & $t_0XX^{-1}$ & $t_0X^{-1}X$ & $t_0YY^{-1}$ & $t_0Y^{-1}Y$
\\
\hline
Resolution & $t_0$ & $t_0$ & $t_0$ & $t_0$
\end{tabular}
}

We conclude that every ambiguity is resolvable, so by the Diamond Lemma~\cite[Theorem~1.2]{berg}
the irreducible words form a~basis for $\hat H_q$.
The result follows.
\end{proof}

In Proposition~\ref{prop:basisv1} we gave a~basis for $\hat H_q$.
In Proposition~\ref{prop:basisvar} below we give a~variation on this basis.

Let $\lambda$ denote an indeterminate.
Let $\mathbb F\lbrack\lambda,\lambda^{-1}\rbrack$ denote the $\mathbb F$-algebra of Laurent polynomials in
$\lambda$ that have all coef\/f\/icients in $\mathbb F$.
\begin{Lemma}
\label{lem:basisfact}
The following is a~basis for the $\mathbb F$-vector space $\mathbb F\lbrack\lambda,\lambda^{-1}\rbrack$:
\begin{gather}
\lambda^k(\lambda+\lambda^{-1})^{\ell},
\qquad
k\in\lbrace0,1\rbrace,
\qquad
\ell\in\mathbb N.
\label{eq:lbasis}
\end{gather}
\end{Lemma}

\begin{proof}
The vectors $\lbrace\lambda^i\rbrace_{i\in\mathbb Z}$ form a~basis for the $\mathbb F$-vector space
$\mathbb F\lbrack\lambda,\lambda^{-1}\rbrack$.
List the elements of this basis in the following order:
\begin{gather}
\label{eq:order}
1,
\quad
\lambda,
\quad
\lambda^{-1},
\quad
\lambda^2,
\quad
\lambda^{-2},
\quad
\lambda^3,
\quad
\lambda^{-3},
\quad
\ldots.
\end{gather}
List the elements of~\eqref{eq:lbasis} in the following order:
\begin{gather}
\label{eq:order1}
1,
\quad
\lambda,
\quad
\lambda+\lambda^{-1},
\quad
\lambda\big(\lambda+\lambda^{-1}\big),
\quad
\big(\lambda+\lambda^{-1}\big)^2,
\quad
\lambda\big(\lambda+\lambda^{-1}\big)^2,
\quad
\ldots.
\end{gather}
Write each element of~\eqref{eq:order1} as a~linear combination of~\eqref{eq:order}.
Consider the corresponding coef\/f\/icient matrix.
This matrix is upper triangular with all diagonal entries~1.
The result follows.
\end{proof}

For a~subset $S$ of any algebra let $\langle S\rangle$ denote the subalgebra generated by~$S$.
\begin{Definition}
\label{def:calx}
Let $\mathbb T$ denote the following subalgebra of $\hat H_q$:
\begin{gather*}
\mathbb T=\langle t^{\pm1}_0,T_1,T_2,T_3\rangle.
\end{gather*}
\end{Definition}

Let $\lbrace\lambda_i\rbrace_{i=0}^3$ denote mutually commuting indeterminates.
By construction the $\mathbb F$-algebra $\mathbb T$ is commutative and generated by
$t^{\pm1}_0$, $T_1$, $T_2$, $T_3$.
Therefore there exists a~surjective $\mathbb F$-algebra homomorphism $\varphi:\mathbb
F\lbrack\lambda^{\pm1}_0,\lambda_1,\lambda_2,\lambda_3\rbrack\to\mathbb T$ that sends
\begin{gather*}
\lambda^{\pm1}_0\mapsto t^{\pm1}_0,
\qquad
\lambda_1\mapsto T_1,
\qquad
\lambda_2\mapsto T_2,
\qquad
\lambda_3\mapsto T_3.
\end{gather*}
\begin{Proposition}\label{prop:2parts}
The above homomorphism $\varphi$ is an isomorphism.
Moreover, in each line below the displayed vectors form a~basis for the $\mathbb F$-vector space~$\mathbb T$:
\begin{gather}\label{eq:basiscalt}
t_0^k T_0^\ell T_1^r T_2^s T_3^t,
\qquad
k\in\lbrace0,1\rbrace,
\qquad
\ell,r,s,t\in\mathbb N;
\\
t_0^k T_1^r T_2^s T_3^t,
\qquad
k\in\mathbb Z,
\qquad
r,s,t\in\mathbb N.
\label{eq:basiscalt2}
\end{gather}
\end{Proposition}

\begin{proof}
By Lemma~\ref{lem:basisfact} the following is a~basis for the $\mathbb F$-vector space $\mathbb
F\lbrack\lambda^{\pm1}_0,\lambda_1,\lambda_2,\lambda_3\rbrack$:
\begin{gather}
\lambda_0^k\big(\lambda_0+\lambda^{-1}_0\big)^\ell\lambda_1^r\lambda_2^s\lambda_3^t,
\qquad
k\in\lbrace0,1\rbrace,
\qquad
\ell,r,s,t\in\mathbb N.
\label{eq:modbase}
\end{gather}
The homomorphism $\varphi$ sends the vectors~\eqref{eq:modbase} to the vectors~\eqref{eq:basiscalt};
therefore the vectors~\eqref{eq:basiscalt} span~$\mathbb T$.
The vectors~\eqref{eq:basiscalt} are linearly independent by Proposition~\ref{prop:basisv1}.
Therefore the vectors~\eqref{eq:basiscalt} form a~basis for~$\mathbb T$.
Consequently $\varphi$ is an isomorphism and~\eqref{eq:basiscalt2} is a~basis for~$\mathbb T$.
\end{proof}

Recall the elements $\alpha$, $\beta$, $\gamma$ of $\hat H_q$
from~\eqref{eq:alphainH}--\eqref{eq:gammainH}.
By those equations $\alpha$, $\beta$, $\gamma$ are contained in $\mathbb T$.
More precisely,~\eqref{eq:alphainH}--\eqref{eq:gammainH} show how $\alpha$, $\beta$, $\gamma$ look in the basis
for $\mathbb T$ from~\eqref{eq:basiscalt2}.
The elements $\alpha$, $\beta$, $\gamma$ look as follows in the basis for $\mathbb T$ from~\eqref{eq:basiscalt}:
\begin{gather*}
\alpha=qT_0T_1-\big(q-q^{-1}\big)t_0T_1+T_2T_3,
\qquad
\beta=qT_0T_3-\big(q-q^{-1}\big)t_0T_3+T_1T_2,
\\
\gamma=qT_0T_2-\big(q-q^{-1}\big)t_0T_2+T_3T_1.
\end{gather*}

We now consider the subalgebras $\langle X^{\pm1}\rangle$ and $\langle Y^{\pm1}\rangle$ of $\hat
H_q$.
By Proposition~\ref{prop:basisv1} the vectors $\lbrace X^i\rbrace_{i\in\mathbb Z}$ form a~basis for
$\langle X^{\pm1}\rangle$ and the vectors $\lbrace Y^i\rbrace_{i\in\mathbb Z}$ form a~basis for $\langle
Y^{\pm1}\rangle$.

\begin{Lemma}
\label{lem:xyiso}
There exists an isomorphism of $\mathbb F$-algebras $\mathbb F\lbrack\lambda^{\pm1}\rbrack\to\langle
X^{\pm1}\rangle$ that sends $\lambda\mapsto X$.
There exists an isomorphism of $\mathbb F$-algebras $\mathbb F\lbrack\lambda^{\pm1}\rbrack\to\langle
Y^{\pm1}\rangle$ that sends $\lambda\mapsto Y$.
\end{Lemma}
\begin{Proposition}
\label{prop:factoriz}
The $\mathbb F$-linear map
\begin{gather*}
\langle Y^{\pm1}\rangle\otimes\langle X^{\pm1}\rangle\otimes\mathbb T\rightarrow\hat H_q,
\qquad
u\otimes v\otimes w
\mapsto uvw
\end{gather*}
is a~bijection.
\end{Proposition}
\begin{proof}
By Proposition~\ref{prop:basisv1}, Lemma~\ref{lem:xyiso}, and since~\eqref{eq:basiscalt} is a~basis for
$\mathbb T$.
\end{proof}

We now give a~variation on the basis for $\hat H_q$ given in Proposition~\ref{prop:basisv1}.
\begin{Proposition}
\label{prop:basisvar}
The following is a~basis for the $\mathbb F$-vector space $\hat H_q$:
\begin{gather}
Y^i X^j t_0^k T_1^r T_2^s T_3^t,
\qquad
i,j,k\in\mathbb Z,
\qquad
r,s,t\in\mathbb N.
\label{eq:basisxy2}
\end{gather}
\end{Proposition}
\begin{proof}
By Proposition~\ref{prop:factoriz} and since~\eqref{eq:basiscalt2} is a~basis for $\mathbb T$.
\end{proof}

\section{The coef\/f\/icient matrix}\label{Section8}

Suppose we have an element of $\hat H_q$ that we wish to express as a~linear combination of the
vectors~\eqref{eq:basisxy1} or~\eqref{eq:basisxy2}.
In order to describe the result ef\/f\/iciently we will use the following notation.
\begin{Definition}
\label{def:coeffmatrix}
  By Proposition~\ref{prop:factoriz} each $h\in\hat H_q$ can be written as
\begin{gather*}
h=\sum_{i,j\in\mathbb Z}Y^iX^j t_{ij},
\qquad
t_{ij}\in\mathbb T.
\end{gather*}
Moreover for $i,j\in\mathbb Z$ the element $t_{ij}$ is uniquely determined by $h$.
We call $t_{ij}$ the {\it coefficient of $Y^iX^j$ in $h$}.
The {\it coefficient matrix} for $h$ has rows and columns indexed by $\mathbb Z$ and $(i,j)$-entry $t_{ij}$
for $i,j\in\mathbb Z$.
We view
\begin{gather*}
h:
\qquad
\begin{tabular}{c|ccccccc}
&$\cdots$&$X^{-2}$&$X^{-1}$&$1$&$X$&$X^2$&$\cdots$
\\
\hline
$\vdots$&&&&$\vdots$&&&
\\
$Y^{-2}$&&$t_{-2,-2}$&$t_{-2,-1}$&$t_{-2,0}$&$t_{-2,1}$&$t_{-2,2}$&
\\
$Y^{-1}$&&$t_{-1,-2}$&$t_{-1,-1}$&$t_{-1,0}$&$t_{-1,1}$&$t_{-1,2}$
\\
$1$&$\cdots$&$t_{0,-2}$&$t_{0,-1}$&$t_{0,0}$&$t_{0,1}$&$t_{0,2}$&$\cdots$
\\
$Y$&&$t_{1,-2}$&$t_{1,-1}$&$t_{1,0}$&$t_{1,1}$&$t_{1,2}$&
\\
$Y^2$&&$t_{2,-2}$&$t_{2,-1}$&$t_{2,0}$&$t_{2,1}$&$t_{2,2}$&
\\
$\vdots$&&&&$\vdots$&&&
\\
\end{tabular}
\end{gather*}
A coef\/f\/icient matrix has f\/initely many nonzero entries.
When we display a~coef\/f\/icient matrix, any row or column not shown has all entries zero.
\end{Definition}

\begin{Example}
The coef\/f\/icient matrix for $A$ is
\begin{gather*}
\begin{tabular}{c|ccc}
&$X^{-1}$&$1$&$X$
\\
\hline
$Y^{-1}$&$0$&$1$&$0$\tsep{1pt}
\\
$1$&$0$&$0$&$0$
\\
$Y$&$0$&$1$&$0$
\\
\end{tabular}
\end{gather*}
The coef\/f\/icient matrix for $B$ is
\begin{gather*}
\begin{tabular}
{c|ccc}&$X^{-1}$&$1$&$X$
\\
\hline
$Y^{-1}$&$0$&$0$&$0$\tsep{1pt}
\\
$1$&$1$&$0$&$1$
\\
$Y$&$0$&$0$&$0$
\\
\end{tabular}
\end{gather*}
\end{Example}

Our next goal is to compute the coef\/f\/icient matrix for $C$.
In order to simplify the computation we initially work with an element $\theta\in\hat H_q$ that is closely
related to~$C$.
\begin{Definition}
\label{def:theta}
Def\/ine $\theta\in\hat H_q$ such that
\begin{gather}\label{eq:th}
q C=\gamma-\theta t^{-1}_0,
\end{gather}
where we recall $\gamma=\big(q^{-1}t_0+qt_0^{-1}\big)T_2+T_1T_3$.
\end{Definition}
\begin{Lemma}\label{lem:thetabasis2}
In the basis \eqref{eq:basisxy2} the element $\theta$ looks as follows:
\begin{gather}
\theta=YX^{-1}t_0-Y^{-1}X t^{-1}_0+Y^{-1}T_3+XT_1+q^{-1}t_0^2T_2.
\label{eq:thetabasis2}
\end{gather}
\end{Lemma}
\begin{proof}
Recall that $C=t_0t_2+(t_0t_2)^{-1}$.
We have $t_0t_2=q^{-1}t^{-1}_3T_1-q^{-1}YX^{-1}$ by Lemma~\ref{lem:tprod1}.
Also $t^{-1}_3=T_3-t_3$ and $t_3=XT_0-Xt_0$.
By these comments
\begin{gather*}
t_0t_2=q^{-1}T_1T_3-q^{-1}XT_0T_1+q^{-1}Xt_0T_1-q^{-1}YX^{-1}.
\end{gather*}
We have $(t_0t_2)^{-1}=qt_3T_1-qXY^{-1}$ by Lemma~\ref{lem:tprod1}.
We mentioned $t_3=XT_0-Xt_0$, and the term $XY^{-1}$ can be evaluated using a~reduction rule from
Proposition~\ref{prop:hgenrel}.
The result follows from these observations along with Def\/inition~\ref{def:theta}.
\end{proof}
\begin{Lemma}\label{lem:thetacoef}
The coefficient matrix for $\theta$ is
\begin{gather*}
\begin{tabular}{c|ccc}
&$X^{-1}$&$1$&$X$
\\
\hline
$Y^{-1}$&$0$&$T_3$&$-t_0^{-1}$\tsep{1pt}
\\
$1$&$0$&$q^{-1}t^{2}_0T_2$&$T_1$
\\
$Y$&$t_0$&$0$&$0$
\end{tabular}
\end{gather*}
\end{Lemma}

\begin{proof}
Use Lemma~\ref{lem:thetabasis2}.
\end{proof}
\begin{Lemma}
\label{lem:ccm}
The coefficent matrix for $C$ is
\begin{gather*}
\begin{tabular}{c|ccc}
&$X^{-1}$&$1$&$X$
\\
\hline
$Y^{-1}$&$0$&$-q^{-1}t_0^{-1}T_3$&$q^{-1}t_0^{-2}$\tsep{1pt}
\\
$1$&$0$&$t_0^{-1}T_2+q^{-1}T_1T_3$&$-q^{-1}t_0^{-1}T_1$
\\
$Y$&$-q^{-1}$&$0$&$0$
\end{tabular}
\end{gather*}
\end{Lemma}

\begin{proof}
Use Def\/inition~\ref{def:theta} and Lemma~\ref{lem:thetacoef}.
\end{proof}
\begin{Lemma}\label{lem:xy2calt}
The coefficient matrix for $XC$ is
\begin{gather*}
\begin{tabular}{c|ccccc}
&$X^{-2}$&$X^{-1}$&$1$&$X$&$X^2$
\\
\hline
$Y^{-2}$&$0$&$0$&$0$&$0$&$0$\tsep{1pt}
\\
$Y^{-1}$&$0$&$q^{-3}t_0T_3$&$-q^{-1}T_3^2-q^{-3}t^2_0-q^{-3}$&$q^{-2}(q^{-1}t_0+qt_0^{-1})T_3$&$-q^{-3}$
\\
$1$&$0$&$-q^{-2}t_0T_2$&$q^{-1}t_0T_1+T_2T_3$&$0$&$0$
\\
$Y$&$0$&$0$&$-q$&$0$&$0$
\\
$Y^2$&$0$&$0$&$0$&$0$&$0$
\end{tabular}
\end{gather*}
\end{Lemma}

\begin{proof}
First f\/ind the coef\/f\/icient matrix for $X\theta$.
To do this, in the equation~\eqref{eq:thetabasis2} multiply each term on the left by $X$ and simplify the
result using the reduction rules from Proposition~\ref{prop:hgenrel}.
This yields the coef\/f\/icient matrix for $X\theta$.
Using this coef\/f\/icient matrix and~\eqref{eq:th}, we routinely obtain the coef\/f\/icient matrix for~$XC$.
\end{proof}

We mention two results for later use.
\begin{Lemma}
\label{lem:xic}
We have
\begin{gather*}
X^{-1}C=q^{-2}C\big(X+X^{-1}\big)-XC-q^{-1}\big(q^2-q^{-2}\big)\big(Y+Y^{-1}\big)+q^{-1}\big(q-q^{-1}\big)\alpha,
\end{gather*}
where we recall $\alpha=\big(q^{-1}t_0+qt^{-1}_0\big)T_1+T_2T_3$.
\end{Lemma}
\begin{proof}
In the f\/irst equation of Lemma~\ref{thm:Qlevel2}, eliminate $A$ using $A=Y+Y^{-1}$ and $B$ using
$B=X+X^{-1}$.
In the resulting equation solve for $X^{-1}C$.
\end{proof}
\begin{Lemma}
\label{lem:domain}
Given $h\in\hat H_q$ and $v\in\mathbb T$ such that $hv=0$.
Then $h=0$ or $v=0$.
\end{Lemma}
\begin{proof}
We assume $v\not=0$ and show $h=0$.
Following Def\/inition~\ref{def:coeffmatrix} write
\begin{gather*}
h=\sum_{i,j\in\mathbb Z}Y^iX^j t_{ij},
\qquad
t_{ij}\in\mathbb T.
\end{gather*}
In this equation we multiply each term on the right by $v$ to obtain
\begin{gather*}
0=\sum_{i,j\in\mathbb Z}Y^iX^j t_{ij}v.
\end{gather*}
Note that $t_{ij}v\in\mathbb T$ for $i,j\in\mathbb Z$.
By this and Proposition~\ref{prop:factoriz} we f\/ind $t_{ij}v=0$ for $i,j\in\mathbb Z$.
The algebra $\mathbb T$ is isomorphic to $\mathbb
F\lbrack\lambda^{\pm1}_0,\lambda_1,\lambda_2,\lambda_3\rbrack$ by Proposition~\ref{prop:2parts}.
The algebra $\mathbb F\lbrack\lambda^{\pm1}_0,\lambda_1,\lambda_2,\lambda_3\rbrack$ is a~domain, so
$\mathbb T$ is a~domain.
By this and since $v\not=0$ we f\/ind $t_{ij}=0$ for all $i,j\in\mathbb Z$.
Therefore $h=0$.
\end{proof}

\section{ The proof of Theorem~\ref{thm:main2}}\label{Section9}

In this section we prove Theorem~\ref{thm:main2}.
Recall the Casimir element $\Omega$ in $\Delta_q$, from Def\/inition~\ref{def:casdelta}.
Let $\Omega'$ denote the element~\eqref{eq:omimage}, so that
\begin{gather*}
\Omega'=\big(q+q^{-1}\big)^2-\big(q^{-1}t_0+qt_0^{-1}\big)^2-T_1^2-T_2^2-T_3^2-\big(q^{-1}t_0+qt_0^{-1}\big)T_1T_2T_3.
\end{gather*}
Theorem~\ref{thm:main2} asserts that $\Omega'$ is the image of $\Omega$ under $\psi$.

\begin{proof}[Proof of Theorem~\ref{thm:main2}] By Def\/inition~\ref{def:casdelta} and
Theorem~\ref{thm:main1} the image of $\Omega$ under $\psi$ is the following element of $\hat H_q$:
\begin{gather}
\label{eq:omegach}
q^{-1}A C B+q^{-2}A^2+q^{-2}B^2+q^{2}C^2-q^{-1}A\alpha-q^{-1}B\beta-q C\gamma,
\end{gather}
where $\alpha$, $\beta$, $\gamma$ are from~\eqref{eq:alphainH}--\eqref{eq:gammainH}.
We show that~\eqref{eq:omegach} is equal to~$\Omega'$.
Def\/ine $D$ to be~\eqref{eq:omegach} minus~$\Omega'$.
We show that $D=0$.
Our strategy is to f\/ind the coef\/f\/icient matrix for $D$ in the sense of
Def\/inition~\ref{def:coeffmatrix}.
Using $A=Y+Y^{-1}$ and $B=X+X^{-1}$ we obtain
\begin{gather}
D=q^{-1}\big(Y+Y^{-1}\big)C\big(X+X^{-1}\big)+q^{-2}\big(Y+Y^{-1}\big)^2+q^{-2}\big(X+X^{-1}\big)^2
\nonumber
\\
\hphantom{D=}{}
+ q^{2}C^2-q^{-1}\big(Y+Y^{-1}\big)\alpha-q^{-1}\big(X+X^{-1}\big)\beta-q C\gamma-\Omega'.
\label{eq:Dform}
\end{gather}
In order to evaluate $D$ further we consider the term $C^2$.
In this product eliminate the f\/irst factor using the formula for $C$ from Lemma~\ref{lem:ccm}.
Simplify the result using the fact that $C$ commutes with~$t_0$; this gives
\begin{gather*}
C^2=q^{-1}Y^{-1}X C t_0^{-2}-q^{-1}YX^{-1}C-q^{-1}X C t_0^{-1}T_1
\\
\hphantom{C^2=}{}
- q^{-1}Y^{-1}C t_0^{-1}T_3+C\big(t_0^{-1}T_2+q^{-1}T_1T_3\big).
\end{gather*}
In the above formula we eliminate $X^{-1}C$ using Lemma~\ref{lem:xic}.
Evaluating~\eqref{eq:Dform} using the results we obtain
\begin{gather}
D=q C\big(T_1T_3-\gamma+qt_0^{-1}T_2\big)-q Y^{-1}C t_0^{-1}T_3+q^{-1}Y^{-1}C\big(X+X^{-1}\big)
\nonumber
\\
\hphantom{D=}{}
- q XC t_0^{-1}T_1+q YXC+q Y^{-1}XC t_0^{-2}+G,
\label{lem:Dformf}
\end{gather}
where
\begin{gather}
G=q^2Y^2+q^{-2}Y^{-2}-q Y\alpha-q^{-1}Y^{-1}\alpha+q^{-2}X^2+q^{-2}X^{-2}
\nonumber
\\
\hphantom{G=}{}
- q^{-1}X\beta-q^{-1}X^{-1}\beta+q^2+3q^{-2}-\Omega'.
\label{eq:Rform}
\end{gather}
We continue to compute the coef\/f\/icient matrix of $D$.
For the next step we will display the coef\/f\/icient matrix for a~number of elements in $\hat H_q$.
When we display these coef\/f\/icient matrices we just display the $(i,j)$ entry for $-2\leq i,j\leq2$,
since it turns out that all the other entries are zero.
Consider the element $C$ of $\hat H_q$.
By Lemma~\ref{lem:ccm} the coef\/f\/icient matrix for $C$ is
\begin{gather}\label{eq:cmat}
\begin{tabular}{c|ccccc}
&$X^{-2}$&$X^{-1}$&$1$&$X$&$X^2$
\\
\hline
$Y^{-2}$&0&0&0&0&0\tsep{1pt}
\\
$Y^{-1}$&0&0&$-q^{-1}t_0^{-1}T_3$&$q^{-1}t_0^{-2}$&0
\\
$1$&0&0&$t_0^{-1}T_2+q^{-1}T_1T_3$&$-q^{-1}t_0^{-1}T_1$&0
\\
$Y$&0&$-q^{-1}$&0&0&0
\\
$Y^2$&0&0&0&0&0
\\
\end{tabular}
\end{gather}

The coef\/f\/icient matrix for $Y^{-1}C$ is
\begin{gather}\label{eq:yicmat}
\begin{tabular}{c|ccccc}
&$X^{-2}$&$X^{-1}$&$1$&$X$&$X^2$
\\
\hline
$Y^{-2}$&0&0&$-q^{-1}t_0^{-1}T_3$&$q^{-1}t_0^{-2}$&0\tsep{1pt}
\\
$Y^{-1}$&0&0&$t_0^{-1}T_2+q^{-1}T_1T_3$&$-q^{-1}t_0^{-1}T_1$&0
\\
$1$&0&$-q^{-1}$&0&0&0
\\
$Y$&0&0&0&0&0
\\
$Y^2$&0&0&0&0&0
\end{tabular}
\end{gather}

By this and since $t_0$ commutes with $X+X^{-1}$, the coef\/f\/icient matrix for
$Y^{-1}C(X+X^{-1})$ is
\begin{gather}\label{eq:yicxpxi}
\begin{tabular}{c|ccccc}
&$X^{-2}$&$X^{-1}$&$1$&$X$&$X^2$
\\
\hline
$Y^{-2}$&0&$-q^{-1}t_0^{-1}T_3$&$q^{-1}t_0^{-2}$&$-q^{-1}t_0^{-1}T_3$&$q^{-1}t_0^{-2}$\tsep{2pt}
\\
$Y^{-1}$&0&$t_0^{-1}T_2+q^{-1}T_1T_3$&$-q^{-1}t_0^{-1}T_1$&$t_0^{-1}T_2+q^{-1}T_1T_3$&$-q^{-1}t_0^{-1}T_1$
\\
$1$&$-q^{-1}$&0&$-q^{-1}$&0&0
\\
$Y$&0&0&0&0&0
\\
$Y^2$&0&0&0&0&0
\end{tabular}
\end{gather}

By Lemma~\ref{lem:xy2calt} the coef\/f\/icient matrix for $XC$ is
\begin{gather}\label{eq:xc}
\begin{tabular}{c|ccccc}
&$X^{-2}$&$X^{-1}$&$1$&$X$&$X^2$
\\
\hline
$Y^{-2}$&0&0&0&0&0\tsep{1pt}
\\
$Y^{-1}$&0&$q^{-3}t_0T_3$&$-q^{-1}T_3^2-q^{-3}t^2_0-q^{-3}$&$q^{-2}\big(q^{-1}t_0+qt_0^{-1}\big)T_3$&$-q^{-3}$
\\
$1$&0&$-q^{-2}t_0T_2$&$q^{-1}t_0T_1+T_2T_3$&0&0
\\
$Y$&0&0&$-q$&0&0
\\
$Y^2$&0&0&0&0&0
\end{tabular}
\end{gather}

The coef\/f\/icient matrix for $YXC$ is
\begin{gather}\label{eq:yxc}
\begin{tabular}{c|ccccc}
&$X^{-2}$&$X^{-1}$&$1$&$X$&$X^2$
\\
\hline
$Y^{-2}$&0&0&0&0&0\tsep{1pt}
\\
$Y^{-1}$&0&0&0&0&0
\\
$1$&0&$q^{-3}t_0T_3$&$-q^{-1}T_3^2-q^{-3}t^2_0-q^{-3}$&$q^{-2}\big(q^{-1}t_0+qt_0^{-1}\big)T_3$&$-q^{-3}$
\\
$Y$&0&$-q^{-2}t_0T_2$&$q^{-1}t_0T_1+T_2T_3$&0&0
\\
$Y^2$&0&0&$-q$&0&0
\end{tabular}
\end{gather}
The coef\/f\/icient matrix for $Y^{-1}XC$ is
\begin{gather}\label{eq:yixc}
\begin{tabular}{c|ccccc}
&$X^{-2}$&$X^{-1}$&$1$&$X$&$X^2$
\\
\hline
$Y^{-2}$&0&$q^{-3}t_0T_3$&$-q^{-1}T_3^2-q^{-3}t^2_0-q^{-3}$&$q^{-2}\big(q^{-1}t_0+qt_0^{-1}\big)T_3$&$-q^{-3}$\tsep{2pt}
\\
$Y^{-1}$&0&$-q^{-2}t_0T_2$&$q^{-1}t_0T_1+T_2T_3$&0&0
\\
$1$&0&0&$-q$&0&0
\\
$Y$&0&0&0&0&0
\\
$Y^2$&0&0&0&0&0
\end{tabular}
\end{gather}
By~\eqref{eq:Rform} the coef\/f\/icient matrix for $G$ is
\begin{gather}\label{eq:rest}
\begin{tabular}{c|ccccc}
&$X^{-2}$&$X^{-1}$&$1$&$X$&$X^2$
\\
\hline
$Y^{-2}$&0&0&$q^{-2}$&0&0\tsep{1pt}
\\
$Y^{-1}$&0&0&$-q^{-1}\alpha$&0&0
\\
$1$&$q^{-2}$&$-q^{-1}\beta$&$q^2+3q^{-2}-\Omega'$&$-q^{-1}\beta$&$q^{-2}$
\\
$Y$&0&0&$-q\alpha$&0&0
\\
$Y^2$&0&0&$q^2$&0&0
\end{tabular}
\end{gather}
We now evaluate~\eqref{lem:Dformf} using~\eqref{eq:cmat}--\eqref{eq:rest}.
One routinely checks that~\eqref{eq:cmat} times $q(T_1T_3-\gamma+qt_0^{-1}T_2)$ minus~\eqref{eq:yicmat}
times $q t_0^{-1}T_3$ plus~\eqref{eq:yicxpxi} times~$q^{-1}$ minus~\eqref{eq:xc} times~$qt_0^{-1}T_1$
plus~\eqref{eq:yxc} times~$q$ plus~\eqref{eq:yixc} times~$q t_0^{-2}$ plus~\eqref{eq:rest} is equal to zero.
Evaluating~\eqref{lem:Dformf} in this light we f\/ind that the coef\/f\/icient matrix of~$D$ is zero.
Therefore $D=0$ and the result follows.
\end{proof}

From now on we retain the notation $\Omega$ for its image under the map $\psi:\Delta_q\to\hat
H_q$.
Thus the element $\Omega$ of $\hat H_q$ satisf\/ies
\begin{gather}
\Omega=
\label{eq:omegadisp2}
\big(q+q^{-1}\big)^2-\big(q^{-1}t_0+qt_0^{-1}\big)^2-T_1^2-T_2^2-T_3^2-\big(q^{-1}t_0+qt_0^{-1}\big)T_1T_2T_3.
\end{gather}

\section{Some results concerning algebraic independence}\label{Section10}

Our next general goal is to prove Theorem~\ref{thm:main3}.
The proof will be completed in Section~\ref{Section12}.
In the present section we establish some results about algebraic independence that will be used in the
proof.

Let $\lbrace x_i\rbrace_{i=1}^4$ denote mutually commuting indeterminates.
Motivated by the form of~\eqref{eq:alphainH}--\eqref{eq:gammainH} and~\eqref{eq:omegadisp2} we consider the
following elements in $\mathbb F\lbrack x_1,x_2,x_3,x_4\rbrack$:
\begin{gather}
y_1=x_1x_2x_3x_4+x_1^2+x_2^2+x_3^2+x_4^2,
\label{eq:y1def}
\\
y_2=x_1x_2+x_3x_4,
\qquad
y_3=x_1x_3+x_2x_4,
\qquad
y_4=x_1x_4+x_2x_3.
\label{eq:y234def}
\end{gather}
\begin{Lemma}[\protect{\cite[Lemma~8.1]{uawe}}]
\label{lem:fourindep}
The elements $\lbrace y_i\rbrace_{i=1}^4$ in \eqref{eq:y1def},~\eqref{eq:y234def} are algebraically
independent over $\mathbb F$.
\end{Lemma}

Recall the algebra $\mathbb T$ from Def\/inition~\ref{def:calx}.
\begin{Lemma}\label{cor:algebin4}
The following are algebraically independent elements of $\mathbb T$:
\begin{gather*}
\Omega,
\quad
\alpha,
\quad
\beta,
\quad
\gamma.
\end{gather*}
\end{Lemma}

\begin{proof}
Recall that $\mathbb T$ is generated by $t^{\pm1}_0$, $T_1$, $T_2$, $T_3$.
By Proposition~\ref{prop:2parts} the following are algebraically independent over $\mathbb F$:
\begin{gather*}
t_0,
\quad
T_1,
\quad
T_2,
\quad
T_3.
\end{gather*}
Therefore the following are algebraically independent over $\mathbb F$:
\begin{gather}
q^{-1}t_0+qt_0^{-1},
\quad
T_1,
\quad
T_2,
\quad
T_3.
\label{eq:Tseq}
\end{gather}
Denote the sequence~\eqref{eq:Tseq} by $\lbrace X_i\rbrace_{i=1}^4$.
By Lemma~\ref{lem:fourindep} the following are algebraically independent over $\mathbb F$:
\begin{gather*}
X_1X_2X_3X_4+X_1^2+X_2^2+X_3^2+X_4^2,
\qquad
X_1X_2+X_3X_4,
\\
X_1X_3+X_2X_4,
\qquad
X_1X_4+X_2X_3.
\end{gather*}
By~\eqref{eq:alphainH}--\eqref{eq:gammainH} and~\eqref{eq:omegadisp2} the above four elements are
\begin{gather*}
\big(q+q^{-1}\big)^2-\Omega,
\quad
\alpha,
\quad
\beta,
\quad
\gamma.
\end{gather*}
The result follows.
\end{proof}
\begin{Definition}
  Let $\mathbb P$ denote the following subalgebra of $\mathbb T$:
\begin{gather*}
\mathbb P=\langle\Omega,\alpha,\beta,\gamma\rangle.
\end{gather*}
\end{Definition}

We set some notation.
For subspaces $U$, $V$ of $\hat H_q$ def\/ine $UV={\rm Span}_\mathbb F\lbrace uv\,|\,u\in U,\,v\in V\rbrace$.

In order to motivate the next few sections let us brief\/ly return to the map
$\psi:\Delta_q\to\hat H_q$ from Theorem~\ref{thm:main1}.
Our current goal is to show that $\psi$ is injective.
Recall that $\Delta_q$ is generated by $A$, $B$, $C$.
Therefore the image of $\Delta_q$ under $\psi$ is the subalgebra $\langle A,B,C\rangle$ of $\hat H_q$.
By Theorem~\ref{thm:usefulbasis} the vectors~\eqref{eq:usefulbasis} form a~basis for $\Delta_q$.
Applying $\psi$ to this basis, we f\/ind that the following vectors span $\langle A,B,C\rangle$:
\begin{gather}
A^i C^j B^k\Omega^\ell\alpha^r\beta^s\gamma^t,
\qquad
j\in\lbrace0,1\rbrace,
\qquad
i,k,\ell,r,s,t\in\mathbb N.
\label{eq:usefulbasis2}
\end{gather}
Consequently
\begin{gather}
\langle A,B,C\rangle=\langle A\rangle\langle B\rangle\mathbb P+\langle A\rangle C\langle B\rangle\mathbb P.
\label{eq:abp}
\end{gather}

In order to show that $\psi$ is injective, it suf\/f\/ices to show that the
vectors~\eqref{eq:usefulbasis2} are linearly independent.
To show this, it will be convenient to expand our focus from the algebra $\langle A,B,C\rangle$ to the
algebra $\langle A,B,C,\mathbb T\rangle=\langle A,B,C,t_0^{\pm1},T_1,T_2,T_3\rangle$.
By~\eqref{eq:abp}, and since everything in $\langle A,B,C\rangle$ commutes with everything in~$\mathbb T$,
\begin{gather}
\langle A,B,C,\mathbb T\rangle=\langle A\rangle\langle B\rangle\mathbb T+\langle A\rangle C\langle B\rangle\mathbb T.
\label{eq:abt}
\end{gather}
We will show that the following is a~basis for the $\mathbb F$-vector space $\langle A,B,C,\mathbb
T\rangle$:
\begin{gather*}
A^iC^jB^k t_0^\ell T_1^rT_2^s T_3^t,
\qquad
j\in\lbrace0,1\rbrace,
\qquad
\ell\in\mathbb Z,
\qquad
i,k,r,s,t\in\mathbb N.
\end{gather*}
It will follow from this and Lemma~\ref{cor:algebin4} that the vectors~\eqref{eq:usefulbasis2} are linearly
independent.

\section[The structure of $\hat H_q$]{The structure of $\boldsymbol{\hat H_q}$}\label{Section11}

In this section we establish some results about $\hat H_q$ that will be used in the proof of
Theorem~\ref{thm:main3}.
Recall $A=Y+Y^{-1}$ and $B=X+X^{-1}$.
\begin{Lemma}
\label{lem:abbasis}
The following is a~basis for $\langle Y^{\pm1}\rangle$:
\begin{gather*}
Y^k A^\ell,
\qquad
k\in\lbrace0,1\rbrace,
\qquad
\ell\in\mathbb N.
\end{gather*}
The following is a~basis for $\langle X^{\pm1}\rangle$:
\begin{gather*}
X^k B^\ell,
\qquad
k\in\lbrace0,1\rbrace,
\qquad
\ell\in\mathbb N.
\end{gather*}
\end{Lemma}
\begin{proof}
Combine Lemma~\ref{lem:basisfact} and Lemma~\ref{lem:xyiso}.
\end{proof}
\begin{Lemma}
\label{lem:xydec}
The following sums are direct:
\begin{gather*}
\langle Y^{\pm1}\rangle=\langle A\rangle+Y\langle A\rangle,
\qquad
\langle X^{\pm1}\rangle=\langle B\rangle+X\langle B\rangle.
\end{gather*}
For each summand a~basis is given in the table below.
\begin{gather*}
\begin{tabular}{c|c}
{\rm subspace}&{\rm basis}
\\
\hline
$\langle A\rangle$&$A^i\tsep{1pt}
\qquad
i\in\mathbb N$
\\
$Y\langle A\rangle$&$YA^i
\qquad
i\in\mathbb N$
\\
$\langle B\rangle$&$B^i
\qquad
i\in\mathbb N$
\\
$X\langle B\rangle$&$XB^i
\qquad
i\in\mathbb N$
\end{tabular}
\end{gather*}
\end{Lemma}

\begin{proof}
Use Lemma~\ref{lem:abbasis}.
\end{proof}
\begin{Proposition}
\label{prop:fourdec}
The following sum is direct:
\begin{gather}
\label{eq:dec1}
\hat H_q=\langle A\rangle\langle B\rangle\mathbb T+\langle A\rangle X\langle B\rangle\mathbb T
+\langle A\rangle Y\langle B\rangle\mathbb T+\langle A\rangle YX\langle B\rangle\mathbb T.
\end{gather}
For each summand a~basis is given in the table below.
\begin{gather*}
\begin{tabular}{c|c}
{\rm subspace}&{\rm basis}
\\
\hline
$\langle A\rangle\langle B\rangle\mathbb T$&$A^i B^j t_0^k T_1^r T_2^s T_3^t
\qquad
k\in\mathbb Z,
\qquad
i,j,r,s,t\in\mathbb N$\tsep{1pt}
\\
$\langle A\rangle X\langle B\rangle\mathbb T$&$A^i X B^j t_0^k T_1^r T_2^s T_3^t
\qquad
k\in\mathbb Z,
\qquad
i,j,r,s,t\in\mathbb N$
\\
$\langle A\rangle Y\langle B\rangle\mathbb T$&$A^i Y B^j t_0^k T_1^r T_2^s T_3^t
\qquad
k\in\mathbb Z,
\qquad
i,j,r,s,t\in\mathbb N$
\\
$\langle A\rangle YX\langle B\rangle\mathbb T$&$A^iYX B^j t_0^k T_1^r T_2^s T_3^t
\qquad
k\in\mathbb Z,
\qquad
i,j,r,s,t\in\mathbb N$
\end{tabular}
\end{gather*}
\end{Proposition}

\begin{proof}
By Proposition~\ref{prop:factoriz}, Lemma~\ref{lem:xydec}, and since~\eqref{eq:basiscalt2} is a~basis for
$\mathbb T$.
\end{proof}

\begin{Proposition}
\label{prop:tensor4}
For $\nu\in\lbrace1,X,Y,YX\rbrace$ the $\mathbb F$-linear map
\begin{gather*}
\langle A\rangle\otimes\langle B\rangle\otimes\mathbb T\rightarrow\langle A\rangle\nu\langle B\rangle\mathbb T,
\qquad
u\otimes v\otimes w
\mapsto u\nu vw
\end{gather*}
is a~bijection.
\end{Proposition}

\begin{proof}
Use the bases displayed in the table of Proposition~\ref{prop:fourdec}.
\end{proof}

Consider the four summands in the decomposition~\eqref{eq:dec1}.
For each summand we now consider the corresponding projection map.
\begin{Definition}
\label{def:projnu}
For $\nu\in\lbrace1,X,Y,YX\rbrace$ def\/ine an $\mathbb F$-linear map $\pi_\nu:\hat H_q\to\hat H_q$
such that $\pi_\nu$ acts as the identity on $\langle A\rangle\nu\langle B\rangle\mathbb T$, and as 0 on the
other three summands in~\eqref{eq:dec1}.
Thus $\pi_\nu$ is the projection from $\hat H_q$ onto $\langle A\rangle\nu\langle B\rangle\mathbb T$.
For $h\in\hat H_q$ we have
\begin{gather}
\pi_\nu(Ah)=A\pi_\nu(h),
\qquad
\pi_\nu(hB)=\pi_\nu(h)B,
\qquad
\pi_\nu(hv)=\pi_\nu(h)v,
\qquad
 \forall\, v\in\mathbb T.
\label{eq:piprop}
\end{gather}
Moreover
\begin{gather*}
h=\pi_1(h)+\pi_X(h)+\pi_Y(h)+\pi_{YX}(h).
\end{gather*}
\end{Definition}

\begin{Lemma}
\label{lem:ccomp}
For $\nu\in\lbrace1,X,Y,YX\rbrace$ the projections $\pi_\nu(A)$, $\pi_\nu(B)$, $\pi_\nu(C)$ are given in
the table below.
\begin{gather*}
\begin{tabular}{c|ccc}
$\nu$&$\pi_\nu(A)$&$\pi_\nu(B)$&$\pi_\nu(C)$
\\
\hline
$1$&$A$&$B$&$q^{-1}\gamma-q^{-2}t_0T_2-q^{-1}At^{-1}_0T_3$\tsep{1pt}
\\
$X$&$0$&$0$&$q^{-1}AXt^{-2}_0-q^{-1}Xt^{-1}_0T_1$
\\
$Y$&$0$&$0$&$q^{-1}Yt^{-1}_0T_3-q^{-1}YB$
\\
$YX$&$0$&$0$&$q^{-1}YX\big(1-t_0^{-2}\big)$
\end{tabular}
\end{gather*}
\end{Lemma}
\begin{proof}
To get $\pi_\nu(A)$ and $\pi_\nu(B)$, note that each of $A$, $B$ is contained in $\langle A\rangle\langle
B\rangle\mathbb T$.
To get~$\pi_\nu(C)$, consider the formula for $C$ from Lemma~\ref{lem:ccm}.
In this formula eliminate $X^{-1}$, $Y^{-1}$ using $X^{-1}=B-X$ and $Y^{-1}=A-Y$.
\end{proof}

\section{The proof of Theorem~\ref{thm:main3}}\label{Section12}

In this section we will prove Theorem~\ref{thm:main3}.
To prepare for the proof, consider the following subspace of $\hat H_q$:
\begin{gather}
\label{eq:tildeh}
\tilde H_q=\langle A\rangle\langle B\rangle\mathbb T+\langle A\rangle X\langle B\rangle\mathbb T
+\langle A\rangle Y\langle B\rangle\mathbb T+\langle A\rangle YX\langle B\rangle\mathbb T\big(1-t^{-2}_0\big).
\end{gather}
\begin{Lemma}
\label{lem:ds}
The sum in \eqref{eq:tildeh} is direct.
\end{Lemma}
\begin{proof}
Observe that $\mathbb T\big(1-t^{-2}_0\big)$ is contained in $\mathbb T$, so $\langle A\rangle YX\langle
B\rangle\mathbb T\big(1-t^{-2}_0\big)$ is contained in $\langle A\rangle YX\langle B\rangle\mathbb T$.
The result follows in view of Proposition~\ref{prop:fourdec}.
\end{proof}

Note that $\mathbb F\lbrack\lambda^{\pm1}\rbrack(1-\lambda^{-2})$ is an ideal in $\mathbb
F\lbrack\lambda^{\pm1}\rbrack$.
\begin{Lemma}
\label{lem:polyideal}
The following sum is direct:
\begin{gather*}
\mathbb F\lbrack\lambda^{\pm1}\rbrack=\mathbb F1+\mathbb F\lambda^{-1}+\mathbb F\lbrack\lambda^{\pm1}
\rbrack\big(1-\lambda^{-2}\big).
\end{gather*}
In other words, the vectors $1$, $\lambda^{-1}$ form a~basis for a~complement of $\mathbb
F\lbrack\lambda^{\pm1}\rbrack(1-\lambda^{-2})$ in $\mathbb F\lbrack\lambda^{\pm1}\rbrack$.
\end{Lemma}

\begin{proof}
One checks that the vectors
\begin{gather*}
1,
\quad
\lambda^{-1},
\quad
1-\lambda^{-2},
\quad
\lambda\big(1-\lambda^{-2}\big),
\quad
\lambda^{-1}(1-\lambda^{-2}),
\quad
\lambda^2(1-\lambda^{-2}),
\quad
\lambda^{-2}(1-\lambda^{-2}),
\quad
\ldots
\end{gather*}
form a~basis for $\mathbb F\lbrack\lambda^{\pm1}\rbrack$.
\end{proof}

Note that $\mathbb T\big(1-t_0^{-2}\big)$ is an ideal in $\mathbb T$.
\begin{Lemma}
\label{lem:basisandc}
The following is a~basis for the $\mathbb F$-vector space $\mathbb T\big(1-t_0^{-2}\big)$:
\begin{gather*}
t_0^k\big(1-t_0^{-2}\big)T^r_1T^s_2T^t_3,
\qquad
k\in\mathbb Z,
\qquad
r,s,t\in\mathbb N.
\end{gather*}
The following is a~basis for a~complement of $\mathbb T\big(1-t_0^{-2}\big)$ in $\mathbb T$:
\begin{gather*}
t^{-k}_0T^r_1T^s_2T^t_3,
\qquad
k\in\lbrace0,1\rbrace,
\qquad
r,s,t\in\mathbb N.
\end{gather*}
\end{Lemma}

\begin{proof}
By Lemma~\ref{lem:polyideal} and the f\/irst assertion of Proposition~\ref{prop:2parts}.
\end{proof}

\begin{Proposition}
\label{prop:basisandc}
The following is a~basis for the $\mathbb F$-vector space $\langle A\rangle YX\langle B\rangle\mathbb
T\big(1-t_0^{-2}\big)$:
\begin{gather}
\label{eq:xytbasis}
A^iYXB^j t_0^k\big(1-t_0^{-2}\big)T^r_1T^s_2T^t_3,
\qquad
k\in\mathbb Z,
\qquad
i,j,r,s,t\in\mathbb N.
\end{gather}
The following is a~basis for a~complement of $\langle A\rangle YX\langle B\rangle\mathbb T\big(1-t_0^{-2}\big)$ in
$\langle A\rangle YX\langle B\rangle\mathbb T$:
\begin{gather*}
A^iYXB^jt^{-k}_0T^r_1T^s_2T^t_3,
\qquad
k\in\lbrace0,1\rbrace,
\qquad
i,j,r,s,t\in\mathbb N.
\end{gather*}
\end{Proposition}
\begin{proof}

Use Proposition~\ref{prop:tensor4} with $\nu=YX$.
Evaluate this using Lemma~\ref{lem:basisandc} along with the fact that $\lbrace A^i\rbrace_{i\in\mathbb N}$
is a~basis for $\langle A\rangle$ and $\lbrace B^i\rbrace_{i\in\mathbb N}$ is a~basis for $\langle
B\rangle$.
\end{proof}
\begin{Corollary}
The following is a~basis for a~complement of $\tilde H_q$ in $\hat H_q$:
\begin{gather*}
A^iYXB^jt^{-k}_0T^r_1T^s_2T^t_3,
\qquad
k\in\lbrace0,1\rbrace,
\qquad
i,j,r,s,t\in\mathbb N.
\end{gather*}
\end{Corollary}
\begin{proof}
This follows from the f\/irst assertion of Proposition~\ref{prop:fourdec}, the def\/inition of $\tilde H_q$
in equation~\eqref{eq:tildeh}, and the last assertion of Proposition~\ref{prop:basisandc}.
\end{proof}
\begin{Lemma}
The following $(i)$--$(iv)$ hold:
\begin{enumerate}\itemsep=0pt
\item[$(i)$] $C\in\tilde H_q$.
\item[$(ii)$] $A\tilde H_q\subseteq\tilde H_q$.
\item[$(iii)$] $\tilde H_q B\subseteq\tilde H_q$.
\item[$(iv)$] $\tilde H_q\mathbb T\subseteq\tilde H_q$.
\end{enumerate}
\end{Lemma}
\begin{proof}
$(i)$ From the column on the right in the table of Lemma~\ref{lem:ccomp}.

$(ii)$, $(iv)$ By equation~\eqref{eq:tildeh}.

$(iii)$ By equation~\eqref{eq:tildeh}, and since $B$ commutes with everything in $\mathbb T$.
\end{proof}

We are about to def\/ine an $\mathbb F$-linear map $\phi:\tilde H_q\to\tilde H_q$.
To def\/ine $\phi$ we give its action on the four summands in~\eqref{eq:tildeh}.
As we will see, the map $\phi$ acts on the f\/irst three summands as a~scalar multiple of the identity.
To give the action of $\phi$ on the fourth summand, we specify what $\phi$ does to the basis for this space
given in~\eqref{eq:xytbasis}.
\begin{Definition}
\label{def:phi}
We def\/ine an $\mathbb F$-linear map $\phi:\tilde H_q\to\tilde H_q$ such that both
\begin{enumerate}\itemsep=0pt
\item[$(i)$] $\phi$ acts as $-q^{-1}$ times the identity on
\begin{gather*}
\langle A\rangle\langle B\rangle\mathbb T+\langle A\rangle X\langle B\rangle\mathbb T+\langle A\rangle Y\langle B\rangle\mathbb T;
\end{gather*}
\item[$(ii)$] for $k\in\mathbb Z$ and $i,j,r,s,t\in\mathbb N$ the map $\phi$ sends
\begin{gather*}
A^iYXB^jt^k_0\big(1-t_0^{-2}\big)T_1^rT_2^s T_3^t\mapsto A^iCB^jt^k_0T_1^rT_2^s T_3^t.
\end{gather*}
\end{enumerate}
\end{Definition}

\begin{Note}
The map $\phi$ is characterized as follows.
Observe that $\phi:\tilde H_q\to\tilde H_q$ is the unique $\mathbb F$-linear map that sends
\begin{gather*}
1\mapsto-q^{-1},
\qquad
X\mapsto-q^{-1}X,
\qquad
Y\mapsto-q^{-1}Y,
\qquad
YX\big(1-t_0^{-2}\big)\mapsto C
\end{gather*}
and satisf\/ies the following for all $h\in\tilde H_q$:
\begin{gather*}
\phi(Ah)=A\phi(h),
\qquad
\phi(hB)=\phi(h)B,
\qquad
\phi(hu)=\phi(h)u,
\qquad
\forall\, u\in\mathbb T.
\end{gather*}
\end{Note}

\begin{Lemma}
We have $\phi^2=q^{-2}1$.
Moreover $\phi$ is a~bijection.
\end{Lemma}
\begin{proof}
The f\/irst assertion is routinely checked using the column on the right in the table of
Lemma~\ref{lem:ccomp}, along with Def\/inition~\ref{def:phi}.
The second assertion is immediate from the f\/irst.
\end{proof}
\begin{Lemma}
\label{lem:images}
Referring to the sum in \eqref{eq:tildeh}, for each summand $U$ the image of $U$ under $\phi$ is
displayed in the table below.
\begin{gather*}
\begin{tabular}{c|c}
$U$&{\rm image of $U$ under $\phi$}
\\
\hline
$\langle A\rangle\langle B\rangle\mathbb T$&$\langle A\rangle\langle B\rangle\mathbb T$
\\
$\langle A\rangle X\langle B\rangle\mathbb T$&$\langle A\rangle X\langle B\rangle\mathbb T$
\\
$\langle A\rangle Y\langle B\rangle\mathbb T$&$\langle A\rangle Y\langle B\rangle\mathbb T$
\\
$\langle A\rangle YX\langle B\rangle\mathbb T\big(1-t^{-2}
_0\big)$&$\langle A\rangle C\langle B\rangle\mathbb T$
\end{tabular}
\end{gather*}
\end{Lemma}

\begin{proof}
Use Def\/inition~\ref{def:phi}.
\end{proof}

\begin{Proposition}
The following sum is direct:
\begin{gather}
\label{eq:tildehalt}
\tilde H_q=\langle A\rangle\langle B\rangle\mathbb T+\langle A\rangle X\langle B\rangle\mathbb T
+\langle A\rangle Y\langle B\rangle\mathbb T+\langle A\rangle C\langle B\rangle\mathbb T.
\end{gather}
Moreover the following is a~basis for the $\mathbb F$-vector space $\langle A\rangle C\langle
B\rangle\mathbb T$:
\begin{gather}
A^iCB^j t_0^k T_1^r T_2^s T_3^t,
\qquad
k\in\mathbb Z,
\qquad
i,j,r,s,t\in\mathbb N.
\label{eq:acbbasis}
\end{gather}
\end{Proposition}

\begin{proof}
The f\/irst assertion is a~consequence of Lemma~\ref{lem:ds} and Lemma~\ref{lem:images}, together with the
fact that $\phi$ is a~bijection.
The second assertion follows from Def\/inition~\ref{def:phi}$(ii)$ and the fact that $\phi$ is a~bijection.
\end{proof}

\begin{Proposition}
\label{prop:close}
The sum \eqref{eq:abt} is direct.
\end{Proposition}
\begin{proof}
The two summands in~\eqref{eq:abt} are included among the four summands in the direct
sum~\eqref{eq:tildehalt}.
\end{proof}

Roughly speaking, the following result amounts to a~universal analog of~\cite[Theorem~2.6]{Koo2}.
\begin{Proposition}
\label{prop:closer}
The following is a~basis for the $\mathbb F$-vector space $\langle A,B,C,\mathbb T\rangle$:
\begin{gather}
A^iC^jB^k t_0^\ell T_1^rT_2^s T_3^t,
\qquad
j\in\lbrace0,1\rbrace,
\qquad
\ell\in\mathbb Z,
\qquad
i,k,r,s,t\in\mathbb N.
\label{eq:finallist}
\end{gather}
\end{Proposition}
\begin{proof}
The set of vectors~\eqref{eq:finallist} consists of the basis for $\langle A\rangle\langle B\rangle\mathbb
T$ from the table of Proposition~\ref{prop:fourdec}, together with the basis for $\langle A\rangle C\langle
B\rangle\mathbb T$ from~\eqref{eq:acbbasis}.
The result follows in view of Proposition~\ref{prop:close}.
\end{proof}

\begin{proof}[Proof of Theorem~\ref{thm:main3}] By Theorem~\ref{thm:usefulbasis} the
vectors~\eqref{eq:usefulbasis} form a~basis for $\Delta_q$.
Applying $\psi$ to this basis, we obtain the following vectors in $\hat H_q$:
\begin{gather*}
A^i C^j B^k\Omega^\ell\alpha^r\beta^s\gamma^t,
\qquad
j\in\lbrace0,1\rbrace,
\qquad
i,k,\ell,r,s,t\in\mathbb N.
\end{gather*}
These vectors are linearly independent by Lemma~\ref{cor:algebin4} and since the
vectors~\eqref{eq:finallist} are linearly independent.
Therefore $\psi$ is injective.
\end{proof}

\section[The elements in $\hat H_q$ that commute with $t_0$]{The elements
in $\boldsymbol{\hat H_q}$ that commute with $\boldsymbol{t_0}$}\label{Section13}

We have now proven the f\/ive theorems from Section~\ref{Section4}.
Recall that these theorems describe the map $\psi:\Delta_q\to\hat H_q$.
Our goal for the remainder of the paper is to obtain three extra results about $\hat H_q$; these results
help to illuminate $\psi$ and may be of independent interest.
The f\/irst extra result concerns the subalgebra $\langle A,B,C,\mathbb T\rangle$ of $\hat H_q$.
This subalgebra was f\/irst mentioned at the end of Section~\ref{Section10}, and a~basis for it was given in
Proposition~\ref{prop:closer}.
Our goal for the present section is to show that
\begin{gather*}
\langle A,B,C,\mathbb T\rangle=\big\lbrace h\in\hat H_q\;|\;t_0h=h t_0\big\rbrace.
\end{gather*}

We will be discussing the $\mathbb F$-linear map $\hat H_q\to\hat H_q$, $h\mapsto t_0h-h
t^{-1}_0$.
\begin{Lemma}
\label{lem:imagedescp}
For $\nu\in\lbrace1,X,Y,YX\rbrace$ the element $t_0\nu-\nu t^{-1}_0$ is given in the table below.
\begin{gather*}
\begin{tabular}
{c|c}
$\nu$&$t_0\nu-\nu t^{-1}_0$
\\
\hline
$1$&$t_0-t^{-1}_0$\tsep{1pt}
\\
$X$&$Bt_0-T_3$
\\
$Y$&$At_0-T_1$
\\
$YX$&$q(Ct_0-T_2)+(AB-T_1T_3)t_0$
\end{tabular}
\end{gather*}
Moreover
\begin{gather}
(AB-T_1T_3)t_0=A(Bt_0-T_3)+(At_0-T_1)t_0T_3-At_0\big(t_0-t^{-1}_0\big)T_3.
\label{eq:ABexpand}
\end{gather}
\end{Lemma}
\begin{proof}
The table is obtained using Lemma~\ref{lem:4rel}.
Equation~\eqref{eq:ABexpand} is routinely checked.
\end{proof}

\begin{Lemma}
\label{lem:imagedesc}
Under the map $h\mapsto t_0h-ht^{-1}_0$ the image of $\hat H_q$ is
\begin{gather*}
\langle A\rangle\big(t_0\!-\!t^{-1}_0\big)\langle B\rangle\mathbb T
+
\langle A\rangle\big(A\!-\!t^{-1}_0T_1\big)\langle B\rangle\mathbb T
+
\langle A\rangle\big(B\!-\!t^{-1}_0T_3\big)\langle B\rangle\mathbb T
+
\langle A\rangle\big(C\!-\!t^{-1}_0T_2\big)\langle B\rangle\mathbb T.
\end{gather*}
This image is contained in $\langle A,B,C,\mathbb T\rangle$.
\end{Lemma}
\begin{proof}
The f\/irst assertion follows from Lemma~\ref{lem:imagedescp}.
The last assertion follows from the f\/irst assertion.
\end{proof}

In~\eqref{eq:tildehalt} we displayed a~direct sum decomposition of $\tilde H_q$.
For each summand we now consider the corresponding projection map.
\begin{Definition}
For $\mu\in\lbrace1,X,Y,C\rbrace$ def\/ine an $\mathbb F$-linear map $P_\mu:\tilde H_q\to\tilde H_q$
such that $P_\mu$ acts as the identity on $\langle A\rangle\mu\langle B\rangle\mathbb T$, and as 0 on the
other three summands in~\eqref{eq:tildehalt}.
Thus $P_\mu$ is the projection from $\tilde H_q$ onto $\langle A\rangle\mu\langle B\rangle\mathbb T$.
For $h\in\tilde H_q$ we have
\begin{gather}
P_\mu(Ah)=A P_\mu(h),
\qquad
P_\mu(hB)=P_\mu(h)B,
\qquad
P_\mu(hv)=P_\mu(h)v,
\qquad
\forall\, v\in\mathbb T.
\label{eq:Pprop}
\end{gather}
Moreover
\begin{gather*}
h=P_1(h)+P_X(h)+P_Y(h)+P_C(h).
\end{gather*}
\end{Definition}

For $h\in\tilde H_q$ we now consider how the projections $P_\mu(h)$ are related to the
projections $\pi_\nu(h)$ from Def\/inition~\ref{def:projnu}.

\begin{Lemma}
\label{lem:pivsP}
Let $h$ denote an element of $\tilde H_q$, and write
\begin{gather}
P_C(h)=\sum_{i,j\in\mathbb N}A^iCB^j t_{ij},
\qquad
t_{ij}\in\mathbb T.
\label{eq:pc}
\end{gather}
Then
\begin{gather*}
\pi_1(h)=P_1(h)+\sum_{i,j\in\mathbb N}A^i\pi_1(C)B^j t_{ij},
\\
\pi_X(h)=P_X(h)+\sum_{i,j\in\mathbb N}A^i\pi_X(C)B^j t_{ij},
\\
\pi_Y(h)=P_Y(h)+\sum_{i,j\in\mathbb N}A^i\pi_Y(C)B^j t_{ij},
\\
\pi_{YX}(h)=q^{-1}\sum_{i,j\in\mathbb N}A^iYX B^j t_{ij}\big(1-t^{-2}_0\big).
\end{gather*}
\end{Lemma}
\begin{proof}
We have both
\begin{gather}
h=\pi_1(h)+\pi_X(h)+\pi_Y(h)+\pi_{YX}(h),
\label{eq:h1}
\\
h=P_1(h)+P_X(h)+P_Y(h)+P_C(h).
\label{eq:h2}
\end{gather}
In~\eqref{eq:h2} eliminate $P_C(h)$ using~\eqref{eq:pc}, and evaluate the result using
\begin{gather*}
C=\pi_1(C)+\pi_X(C)+\pi_Y(C)+\pi_{YX}(C).
\end{gather*}
By Lemma~\ref{lem:ccomp} we have $\pi_{YX}(C)=q^{-1}YX\big(1-t^{-2}_0\big)$.
Subtracting~\eqref{eq:h2} from~\eqref{eq:h1} and using the above comments, we obtain
\begin{gather}
0=\pi_1(h)-P_1(h)-\sum_{i,j\in\mathbb N}A^i\pi_1(C)B^j t_{ij}
\label{eq:t1}
\\
\hphantom{0=}{} + \pi_X(h)-P_X(h)-\sum_{i,j\in\mathbb N}A^i\pi_X(C)B^j t_{ij}
\label{eq:tX}
\\
\hphantom{0=}{} +\pi_Y(h)-P_Y(h)-\sum_{i,j\in\mathbb N}A^i\pi_Y(C)B^j t_{ij}
\label{eq:tY}
\\
\hphantom{0=}{} + \pi_{YX}(h)-q^{-1}\sum_{i,j\in\mathbb N}A^iYX B^j t_{ij}\big(1-t^{-2}_0\big).
\label{eq:tYX}
\end{gather}
The elements~\eqref{eq:t1},~\eqref{eq:tX},~\eqref{eq:tY},~\eqref{eq:tYX} are contained in the subspaces
\begin{gather*}
\langle A\rangle\langle B\rangle\mathbb T,
\qquad
\langle A\rangle X\langle B\rangle\mathbb T,
\qquad
\langle A\rangle Y\langle B\rangle\mathbb T,
\qquad
\langle A\rangle YX\langle B\rangle\mathbb T
\end{gather*}
respectively.
The sum of these subspaces is direct, so each
of~\eqref{eq:t1},~\eqref{eq:tX},~\eqref{eq:tY},~\eqref{eq:tYX} is zero.
The result follows.
\end{proof}
\begin{Lemma}
\label{lem:tiff}
For $h\in\hat H_q$ the following are equivalent:
\begin{enumerate}\itemsep=0pt
\item[$(i)$] $h\in\langle A,B,C,\mathbb T\rangle$.
\item[$(ii)$] $h\big(t_0-t^{-1}_0\big)\in\langle A,B,C,\mathbb T\rangle$.
\end{enumerate}
\end{Lemma}
\begin{proof}
$(i)\Rightarrow (ii)$ Since $t_0-t_0^{-1}\in\mathbb T$.

$(ii)\Rightarrow (i)$ Observe by~\eqref{eq:tildeh} that $h\big(t_0-t^{-1}_0\big)\in\tilde H_q$.
Write
\begin{gather*}
P_C\big(h\big(t_0-t^{-1}_0\big)\big)=\sum_{i,j\in\mathbb N}A^iC B^j t_{ij},
\qquad
t_{ij}\in\mathbb T.
\end{gather*}
We f\/irst show that $h\in\tilde H_q$.
Comparing~\eqref{eq:abt} and~\eqref{eq:tildehalt} we f\/ind
\begin{gather}
P_X\big(h\big(t_0-t^{-1}_0\big)\big)=0,
\qquad
P_Y\big(h\big(t_0-t^{-1}_0\big)\big)=0.
\label{eq:xyzero}
\end{gather}
By the equation on the right in~\eqref{eq:piprop},
\begin{gather*}
\pi_\nu\big(h\big(t_0-t^{-1}_0\big)\big)=\pi_\nu(h)\big(t_0-t_0^{-1}\big),
\qquad
\nu\in\lbrace1,X,Y,YX\rbrace.
\end{gather*}
By this and Lemma~\ref{lem:pivsP},
\begin{gather*}
\pi_{YX}(h)\big(t_0-t^{-1}_0\big)=\pi_{YX}\big(h\big(t_0-t_0^{-1}\big)\big)
=q^{-1}\sum_{i,j\in\mathbb N}A^iYXB^jt_{ij}\big(1-t_0^{-2}\big).
\end{gather*}
By this and Lemma~\ref{lem:domain},
\begin{gather}
\pi_{YX}(h)=q^{-1}\sum_{i,j\in\mathbb N}A^iYXB^jt_{ij}t^{-1}_0.
\label{eq:yxproj}
\end{gather}
In order to show that $h\in\tilde H_q$ we show that $t_0-t^{-1}_0$ divides $t_{ij}$ for all $i,j\in\mathbb N$.
Observe
\begin{gather*}
\pi_Y(h)\big(t_0-t^{-1}_0\big)=\pi_Y\big(h\big(t_0-t_0^{-1}\big)\big)
\overset{\text{by Lemma~\ref{lem:pivsP}}}{=} \; \sum_{i,j\in\mathbb N}A^i\pi_Y(C)B^j t_{ij}
\\
\hphantom{\pi_Y(h)\big(t_0-t^{-1}_0\big)}{}
\overset{\text{by Lemma~\ref{lem:ccomp}}}{=} \;
q^{-1}\sum_{i,j\in\mathbb N}A^i\big(Yt^{-1}_0T_3-YB\big)B^j t_{ij}\\
\hphantom{\pi_Y(h)\big(t_0-t^{-1}_0\big)}{}
=
q^{-1}\sum_{r,s\in\mathbb N}A^r YB^s\big(t_{rs}t^{-1}_0T_3-t_{r,s-1}\big),
\end{gather*}
where $t_{r,-1}=0$ for $r\in\mathbb N$.
From this we see that $t_0-t^{-1}_0$ divides $t_{rs}t^{-1}_0T_3-t_{r,s-1}$ for all $r,s\in\mathbb N$.
By this and induction on $s$ we f\/ind $t_0-t^{-1}_0$ divides $t_{rs}$ for all $r,s\in\mathbb N$.
In other words, for all $r,s\in\mathbb N$ there exists $t'_{rs}\in\mathbb T$ such that
$t_{rs}=t'_{rs}\big(t_0-t_0^{-1}\big)$.
Now using~\eqref{eq:yxproj},
\begin{gather*}
\pi_{YX}(h)=q^{-1}\sum_{i,j\in\mathbb N}A^iYX B^j t'_{ij}\big(1-t_0^{-2}\big)
\in
\langle A\rangle YX\langle B\rangle\mathbb T\big(1-t^{-2}_0\big).
\end{gather*}
By this and~\eqref{eq:tildeh} we f\/ind $h\in\tilde H_q$.
By the equation on the right in~\eqref{eq:Pprop} and the equation on the left in~\eqref{eq:xyzero},
\begin{gather*}
P_X(h)\big(t_0-t^{-1}_0\big)=P_X\big(h\big(t_0-t_0^{-1}\big)\big)=0.
\end{gather*}
Therefore $P_X(h)=0$ in view of Lemma~\ref{lem:domain}.
Similarly
\begin{gather*}
P_Y(h)\big(t_0-t^{-1}_0\big)=P_Y\big(h\big(t_0-t_0^{-1}\big)\big)=0,
\end{gather*}
so $P_Y(h)=0$.
Now
\begin{gather*}
h=P_1(h)+P_X(h)+P_Y(h)+P_C(h)
=P_1(h)+P_C(h)
\\
\hphantom{h=}{}
\in
\langle A\rangle\langle B\rangle\mathbb T+\langle A\rangle C\langle B\rangle\mathbb T
=\langle A,B,C,\mathbb T\rangle.
\end{gather*}
The result follows.
\end{proof}

Roughly speaking, the following result amounts to a~universal analog of~\cite[Theorem~5.1]{Koo2}.
\begin{Theorem}\label{thm:com}
We have
\begin{gather}
\langle A,B,C,\mathbb T\rangle=\big\lbrace h\in\hat H_q\;|\;t_0h=h t_0\big\rbrace.
\label{eq:comt0main}
\end{gather}
\end{Theorem}

\begin{proof}
In~\eqref{eq:comt0main} the inclusion $\subseteq$ holds by Lemma~\ref{cor:txyz} and since $\lbrace
T_i\rbrace_{i=1}^3$ are central in $\hat H_q$.
We now obtain the inclusion $\supseteq$.
Pick $h\in\hat H_q$ such that $t_0h=ht_0$.
We show $h\in\langle A,B,C,\mathbb T\rangle$.
By assumption $t_0h=ht_0$ so $h\big(t_0-t_0^{-1}\big)=t_0h-h t^{-1}_0$.
By Lemma~\ref{lem:imagedesc} $t_0h-h t^{-1}_0\in\langle A,B,C,\mathbb T\rangle$.
By these comments $h\big(t_0-t_0^{-1}\big)\in\langle A,B,C,\mathbb T\rangle$.
By this and Lemma~\ref{lem:tiff} $h\in\langle A,B,C,\mathbb T\rangle$.
\end{proof}

\section[A presentation for the algebra $\langle A,B,C,\mathbb T\rangle$]{A presentation
for the algebra $\boldsymbol{\langle A,B,C,\mathbb T\rangle}$}\label{Section14}

We continue to discuss the subalgebra $\langle A,B,C,\mathbb T\rangle$ of $\hat H_q$.
In this section we give a~presentation for $\langle A,B,C,\mathbb T\rangle$ by generators and relations.
Roughly speaking, this presentation amounts to a~$q$-analog of~\cite[Theorem 2.1]{oblomkov} and a~universal
analog of~\cite[Def\/inition~6.1, Corollary~6.3]{Koo1}.
\begin{Theorem}
\label{thm:t0com}
The $\mathbb F$-algebra $\langle A,B,C,\mathbb T\rangle$ is presented by generators and relations in the
following way.
The generators are $A$, $B$, $C$, $t^{\pm1}_0$, $\lbrace T_i\rbrace_{i=1}^3$.
The relations assert that each of $t^{\pm1}_0$, $\lbrace T_i\rbrace_{i=1}^3$ is central and
$t_0t^{-1}_0=1$, $t^{-1}_0t_0=1$,
\begin{gather*}
A+\frac{qBC-q^{-1}CB}{q^2-q^{-2}}=\frac{\alpha}{q+q^{-1}},
\\
B+\frac{qCA-q^{-1}AC}{q^2-q^{-2}}=\frac{\beta}{q+q^{-1}},
\\
C+\frac{qAB-q^{-1}BA}{q^2-q^{-2}}=\frac{\gamma}{q+q^{-1}},
\\
q^{-1}ACB+q^{-2}A^2+q^{-2}B^2+q^2C^2-q^{-1}A\alpha-\;q^{-1}B\beta-qC\gamma
\\
\qquad{}
=\big(q+q^{-1}\big)^2-\big(q^{-1}t_0+qt^{-1}_0\big)^2-T_1^2-T_2^2-T_3^2-\big(q^{-1}t_0+qt^{-1}_0\big)T_1T_2T_3,
\end{gather*}
where
\begin{gather*}
\alpha=\big(q^{-1}t_0+qt^{-1}_0\big)T_1+T_2T_3,
\qquad
\beta=\big(q^{-1}t_0+qt^{-1}_0\big)T_3+T_1T_2,
\\
\gamma=\big(q^{-1}t_0+qt^{-1}_0\big)T_2+T_3T_1.
\end{gather*}
\end{Theorem}

\begin{proof}
Let $\mathcal A_q$ denote the $\mathbb F$-algebra def\/ined by generators $A$, $B$, $C$, $t^{\pm1}_0$,
$\lbrace T_i\rbrace_{i=1}^3$ and the above relations.
Since these relations hold in $\hat H_q$ there exists an $\mathbb F$-algebra homomorphism $\mathcal
A_q\to\hat H_q$ that sends each generator $A$, $B$, $C$, $t^{\pm1}_0$, $\lbrace T_i\rbrace_{i=1}^3$ of
$\mathcal A_q$ to the corresponding element in $\hat H_q$.
Under this homomorphism the image of $\mathcal A_q$ is the subalgebra $\langle A,B,C,\mathbb T\rangle$ of
$\hat H_q$.
We show that the homomorphism is injective.
To this end, we claim that the following vectors span the $\mathbb F$-vector space $\mathcal A_q$:
\begin{gather}
A^iC^jB^k t_0^\ell T_1^rT_2^s T_3^t,
\qquad
j\in\lbrace0,1\rbrace,
\qquad
\ell\in\mathbb Z,
\qquad
i,k,r,s,t\in\mathbb N.
\label{eq:calaspan}
\end{gather}
To prove the claim, note that the elements $A$, $B$, $C$ of $\mathcal A_q$ satisfy the def\/ining relations for~$\Delta_q$ given in Def\/inition~\ref{def:uaw}.
Therefore there exists an $\mathbb F$-algebra homomorphism $\Delta_q\to\mathcal A_q$ that sends each
generator $A$, $B$, $C$ of $\Delta_q$ to the corresponding element in $\mathcal A_q$.
In~\eqref{eq:usefulbasis} we displayed a~basis for the $\mathbb F$-vector space $\Delta_q$.
When our homomorphism $\Delta_q\to\mathcal A_q$ is applied to a~vector in this basis, the image is
contained in the span of~\eqref{eq:calaspan}.
Therefore the span of~\eqref{eq:calaspan} contains the subalgebra of $\mathcal A_q$ generated by $A$, $B$, $C$.
By construction $\mathcal A_q$ is generated by $A$, $B$, $C$, $t^{\pm1}_0$, $\lbrace T_i\rbrace_{i=1}^3$.
By def\/inition each element $A$, $B$, $C$ of $\mathcal A_q$ commutes with each element $t^{\pm1}_0$, $\lbrace
T_i\rbrace_{i=1}^3$ of $\mathcal A_q$.
By construction the span of~\eqref{eq:calaspan} is closed under multiplication by each element
$t^{\pm1}_0$, $\lbrace T_i\rbrace_{i=1}^3$ of $\mathcal A_q$.
By these comments the vectors~\eqref{eq:calaspan} span $\mathcal A_q$.
The claim is proven.
When we apply our homomorphism $\mathcal A_q\to\hat H_q$ to the vectors~\eqref{eq:calaspan}, we get the
basis for $\langle A,B,C,\mathbb T\rangle$ given in Proposition~\ref{prop:closer}.
Therefore the vectors~\eqref{eq:calaspan} form a~basis for $\mathcal A_q$ and our homomorphism $\mathcal
A_q\to\hat H_q$ is injective.
The result follows.
\end{proof}

\section[The center of $\hat H_q$]{The center of $\boldsymbol{\hat H_q}$}\label{Section15}

In this section we describe the center $Z(\hat H_q)$.

Recall that the $\lbrace T_i\rbrace_{i\in\mathbb I}$ are central in $\hat H_q$.
We are going to show that $\lbrace T_i\rbrace_{i\in\mathbb I}$ generate $Z(\hat H_q)$, provided that $q$ is
not a~root of unity.
In this derivation we will repeatedly use the basis for $\hat H_q$ given in Proposition~\ref{prop:basisv1}.
\begin{Definition}
Let $K$ denote the 2-sided ideal of $\hat H_q$ generated by $\lbrace T_i\rbrace_{i\in\mathbb I}$.
Thus
\begin{gather*}
K=\sum_{i\in\mathbb I}\hat H_q T_i.
\end{gather*}
\end{Definition}

\begin{Lemma}
\label{lem:basisK}
The following is a~basis for the $\mathbb F$-vector space $K$:
\begin{gather*}
Y^iX^j t_0^k T_0^\ell T_1^r T_2^s T_3^t,
\qquad\!
i,j\in\mathbb Z,
\qquad\!
k\in\lbrace0,1\rbrace,
\qquad\!
\ell,r,s,t\in\mathbb N,
\qquad\!
(\ell,r,s,t)\not=(0,0,0,0).
\end{gather*}
\end{Lemma}
\begin{proof}
Use Proposition~\ref{prop:basisv1}.
\end{proof}
\begin{Lemma}
\label{lem:basisCK}
The following is a~basis for a~complement of $K$ in $\hat H_q$:
\begin{gather*}
Y^iX^j t_0^k,
\qquad
i,j\in\mathbb Z,
\qquad
k\in\lbrace0,1\rbrace.
\end{gather*}
\end{Lemma}

\begin{proof}
Compare Proposition~\ref{prop:basisv1} and Lemma~\ref{lem:basisK}.
\end{proof}

\begin{Definition}
\label{def:Kquot}
Let $\overline H_q$ denote the quotient $\mathbb F$-algebra $\overline H_q=\hat H_q/K$.
Recall that the canonical map $\hat H_q\to\overline H_q$ is a~surjective $\mathbb F$-algebra homomorphism
with kernel~$K$.
For $h\in\hat H_q$ let $\overline h$ denote the image of $h$ under this map.
By construction $\overline T_i=0$ for $i\in\mathbb I$.
\end{Definition}
\begin{Lemma}
\label{lem:barbasis}
The following is a~basis for the $\mathbb F$-vector space $\overline H_q$:
\begin{gather*}
\overline Y^i\overline X^j\overline t_0^k,
\qquad
i,j\in\mathbb Z,
\qquad
k\in\lbrace0,1\rbrace.
\end{gather*}
\end{Lemma}

\begin{proof}
Use Lemma~\ref{lem:basisCK}.
\end{proof}

\begin{Lemma}\label{eq:Cbar}
Referring to Definition~{\rm \ref{def:ci}} we have $\overline C_i=0$ for $i\in\mathbb I$.
\end{Lemma}

\begin{proof}
By Proposition~\ref{prop:Ciform} and since $\overline T_j=0$ for $j\in\mathbb I$.
\end{proof}

\begin{Lemma}
\label{lem:overlineh}
The following relations hold in $\overline H_q$:
\begin{gather}
\overline X\,  \overline Y=q^2\overline Y \, \overline X,
\qquad
\overline t_0^2=-1,
\label{eq:rel14}
\\
\overline t_0\overline X=\overline X^{-1}\overline t_0,
\qquad
 \overline t_0\overline Y=\overline Y^{-1}\overline t_0.
\label{eq:rel23}
\end{gather}
\end{Lemma}

\begin{proof}
The equation on the left in~\eqref{eq:rel14} follows from Def\/inition~\ref{def:ci} and Lemma~\ref{eq:Cbar}.
To get the equation on the right in~\eqref{eq:rel14}, apply the map $h\mapsto\overline h$ to each side of
$t_0^2=t_0T_0-1$.
To get the equations in~\eqref{eq:rel23}, apply the map $h\mapsto\overline h$ to each side
of~\eqref{eq:4rel1} and~\eqref{eq:4rel3}.
\end{proof}
\begin{Definition}
\label{def:poset}
We endow the set $\mathbb N^4$ with a~partial order $\leq$ as follows.
Let $(n_0,n_1,n_2,n_3)$ and $(n'_0,n'_1,n'_2,n'_3)$ denote elements of $\mathbb N^4$.
Then $(n_0,n_1,n_2,n_3)\leq(n'_0,n'_1,n'_2,n'_3)$ whenever $n_i\leq n'_i$ for $0\leq i\leq3$.
\end{Definition}

We have some comments.
Fix $(\ell,r,s,t)\in\mathbb N^4$ and def\/ine
\begin{gather}
L=\hat H_q T_0^\ell T_1^r T_2^s T_3^t.
\label{eq:Ldef}
\end{gather}
Then $L$ is a~2-sided ideal of $\hat H_q$ with basis
\begin{gather*}
Y^i X^j t_0^k T_0^{\ell'}T_1^{r'}T_2^{s'}T_3^{t'},
\qquad
i,j\in\mathbb Z,
\qquad
k\in\lbrace0,1\rbrace,
\\
(\ell',r',s',t')\in\mathbb N^4,
\qquad
(\ell,r,s,t)\leq(\ell',r',s',t').
\end{gather*}
Observe that $KL$ is a~2-sided ideal of $\hat H_q$ with basis
\begin{gather*}
Y^i X^j t_0^k T_0^{\ell'}T_1^{r'}T_2^{s'}T_3^{t'},
\qquad
i,j\in\mathbb Z,
\qquad
k\in\lbrace0,1\rbrace,
\\
(\ell',r',s',t')\in\mathbb N^4,
\qquad
(\ell,r,s,t)<(\ell',r',s',t').
\end{gather*}
Def\/ine
\begin{gather}
R=\hat H_q T_0^{\ell+1}+\hat H_q T_1^{r+1}+\hat H_q T_2^{s+1}+\hat H_q T_3^{t+1}.
\label{eq:Rdef}
\end{gather}
Then $R$ is a~2-sided ideal of $\hat H_q$ with basis
\begin{gather*}
Y^iX^j t_0^k T_0^{\ell'}T_1^{r'}T_2^{s'}T_3^{t'},
\qquad
i,j\in\mathbb Z,
\qquad
k\in\lbrace0,1\rbrace,
\\
(\ell',r',s',t')\in\mathbb N^4,
\qquad
(\ell',r',s',t')\not\leq(\ell,r,s,t).
\end{gather*}
Comparing the above bases we f\/ind
\begin{gather}
L\cap R=KL.
\label{lem:GL}
\end{gather}
\begin{Theorem}
\label{thm:centerH}
Assume that $q$ is not a~root of unity.
Then the $\mathbb F$-algebra $Z(\hat H_q)$ is generated by $\lbrace T_i\rbrace_{i\in\mathbb I}$.
\end{Theorem}
\begin{proof}
Consider the subalgebra $\langle T_0,T_1,T_2,T_3\rangle$ of $\hat H_q$.
This subalgebra is contained in $Z(\hat H_q)$.
We assume that the containment is proper, and obtain a~contradiction.
Pick
\begin{gather}
\label{eq:hchoice}
h\in Z\big(\hat H_q\big),
\qquad
h\not\in\langle T_0,T_1,T_2,T_3\rangle.
\end{gather}
In view of Proposition~\ref{prop:basisv1} we write
\begin{gather*}
h=\sum_{\ell,r,s,t\in\mathbb N}h_{\ell,r,s,t}T_0^\ell T_1^r T_2^s T_3^t,
\qquad
h_{\ell,r,s,t}\in\langle Y^{\pm1}\rangle\langle X^{\pm1}\rangle+\langle Y^{\pm1}\rangle\langle X^{\pm1}
\rangle t_0.
\end{gather*}
Def\/ine the set
\begin{gather*}
S(h)=\big\lbrace(\ell,r,s,t)\in\mathbb N^4\;|\;h_{\ell,r,s,t}\not=0\big\rbrace.
\end{gather*}
By construction the cardinality $|S(h)|$ is f\/inite.
Without loss of generality, we assume that $h$ has been chosen such that $|S(h)|$ is minimal subject
to~\eqref{eq:hchoice}.
Note that $h\not=0$ so $S(h)$ is nonempty.
There exists an element of $S(h)$ that is not greater than any other element of $S(h)$, with respect to the
partial order $\leq$ from Def\/inition~\ref{def:poset}.
Denote this element by $(\ell,r,s,t)$.
We will be discussing the corresponding ideals $L$, $R$ of $\hat H_q$ from~\eqref{eq:Ldef}
and~\eqref{eq:Rdef}.
By construction
\begin{gather}
\label{eq:fps}
h-h_{\ell,r,s,t}T_0^\ell T_1^r T_2^s T_3^t\in R.
\end{gather}
Write
\begin{gather}
h_{\ell,r,s,t}=\sum_{i,j\in\mathbb Z}\alpha_{ij}Y^iX^j+\sum_{i,j\in\mathbb Z}\beta_{ij}Y^iX^jt_0,
\qquad
\alpha_{ij}, \beta_{ij}\in\mathbb F.
\label{eq:alphabeta}
\end{gather}
We take the commutator of~\eqref{eq:fps} with each of $X$, $Y$.
We start with $X$.
The ideal $R$ contains
\begin{gather*}
Xh-hX-(Xh_{\ell,r,s,t}-h_{\ell,r,s,t}X)T_0^\ell T_1^r T_2^s T_3^t.
\end{gather*}
By assumption $h\in Z(\hat H_q)$ so $Xh-hX=0$.
Therefore $R$ contains
\begin{gather}
(Xh_{\ell,r,s,t}-h_{\ell,r,s,t}X)T_0^\ell T_1^r T_2^s T_3^t.
\label{eq:hcom}
\end{gather}
The element~\eqref{eq:hcom} is contained in $L$ by~\eqref{eq:Ldef}.
By these comments and~\eqref{lem:GL}, the element~\eqref{eq:hcom} is contained in $KL$.
By Lemma~\ref{lem:domain} the map $\hat H_q\to L$, $g\mapsto gT^\ell_0T^r_1T^s_2T^t_3$ is a~bijection.
Under this map the image of $K$ is $KL$.
Therefore in~\eqref{eq:hcom}, the expression in parenthesis is contained in $K$.
In other words, in the notation of Def\/inition~\ref{def:Kquot},
\begin{gather}
\overline X\,\overline{h}_{\ell,r,s,t}-\overline{h}_{\ell,r,s,t}\,\overline X=0.
\label{eq:comh}
\end{gather}
Expanding~\eqref{eq:comh} using~\eqref{eq:alphabeta} we obtain
\begin{gather*}
0=\sum_{i,j\in\mathbb Z}\alpha_{ij}
\big(\overline X\,\overline Y^i\overline X^j-\overline Y^i\overline X^j\overline X\big)
+\sum_{i,j\in\mathbb Z}\beta_{ij}\big(\overline X\,\overline Y^i\overline X^j\overline t_0
-\overline Y^i\,\overline X^j\overline t_0\,\overline X\big).
\end{gather*}
Simplifying this using Lemma~\ref{lem:overlineh} we obtain
\begin{gather*}
0=\sum_{i,j\in\mathbb Z}\alpha_{ij}\overline Y^i\overline X^{j+1}\big(q^{2i}-1\big)+\sum_{i,j\in\mathbb Z}\beta_{ij}
\big(\overline Y^i\overline X^{j+1}\overline t_0q^{2i}-\overline Y^i\overline X^{j-1}\overline t_0\big).
\end{gather*}
Adjusting the indices $i$, $j$ in the above sums,
\begin{gather*}
0=\sum_{i,j\in\mathbb Z}\overline Y^i\overline X^j\alpha_{i,j-1}\big(q^{2i}-1\big)+\sum_{i,j\in\mathbb Z}
\overline Y^i\overline X^j\overline t_0\big(\beta_{i,j-1}q^{2i}-\beta_{i,j+1}\big).
\end{gather*}
By this and Lemma~\ref{lem:barbasis} we f\/ind
\begin{gather}
\alpha_{i,j-1}\big(q^{2i}-1\big)=0,
\qquad
i,j\in\mathbb Z,
\label{eq:one}
\\
\beta_{i,j-1}q^{2i}-\beta_{i,j+1}=0,
\qquad
i,j\in\mathbb Z.
\label{eq:two}
\end{gather}
Taking the commutator of~\eqref{eq:fps} with $Y$, we similarly obtain
\begin{gather}
\alpha_{i-1,j}\big(q^{2j}-1\big)=0,
\qquad
i,j\in\mathbb Z,
\label{eq:three}
\\
\beta_{i-1,j}-\beta_{i+1,j}q^{-2j}=0,
\qquad
i,j\in\mathbb Z.
\label{eq:four}
\end{gather}
By~\eqref{eq:one},~\eqref{eq:three} and since $q$ is not a~root of unity,
\begin{gather*}
\alpha_{ij}=0
\qquad
\mbox{if}
\quad
(i,j)\not=(0,0),
\quad
i,j\in\mathbb Z.
\end{gather*}
By~\eqref{eq:two} or~\eqref{eq:four}, and since f\/initely many of the~$\beta_{ij}$ are nonzero,
\begin{gather*}
\beta_{ij}=0,
\qquad
i,j\in\mathbb Z.
\end{gather*}
Evaluating~\eqref{eq:alphabeta} using these comments we obtain $h_{\ell,r,s,t}=\alpha_{00}\in\mathbb F$.
Def\/ine
\begin{gather*}
h'=h-h_{\ell,r,s,t}T_0^\ell T_1^r T_2^s T_3^t.
\end{gather*}
We have two comments about $h'$.
First of all,
\begin{gather*}
h-h'\in\langle T_0,T_1,T_2,T_3\rangle\subseteq Z(\hat H_q),
\end{gather*}
so
\begin{gather*}
h'\in Z(\hat H_q),
\qquad
h'\not\in\langle T_0,T_1,T_2,T_3\rangle.
\end{gather*}
Second of all, $S(h')$ is obtained from $S(h)$ by deleting the element $(\ell,r,s,t)$; therefore
$|S(h')|=|S(h)|-1$.
These two comments contradict the minimality of $|S(h)|$.
The result follows.
\end{proof}

Roughly speaking, the following two corollaries amount to a~universal analog
of~\cite[Theorem~5.3]{Koo2}.
\begin{Corollary}
\label{cor:zbasis}
Assume that $q$ is not a~root of unity.
Then the following is a~basis for the $\mathbb F$-vector space $Z(\hat H_q)$:
\begin{gather}
T_0^\ell T_1^r T_2^s T_3^t,
\qquad
\ell,r,s,t\in\mathbb N.
\label{eq:bb}
\end{gather}
\end{Corollary}

\begin{proof}
The vectors~\eqref{eq:bb} span $Z(\hat H_q)$ by Theorem~\ref{thm:centerH}.
The vectors~\eqref{eq:bb} are linearly independent because they are included in the linearly independent
set~\eqref{eq:basiscalt}.
\end{proof}
\begin{Corollary}
Assume that $q$ is not a~root of unity.
Then there exists an isomorphism of $\mathbb F$-algebras $Z(\hat H_q)\to\mathbb
F\lbrack\lambda_0,\lambda_1,\lambda_2,\lambda_3\rbrack$ that sends $T_i\mapsto\lambda_i$ for $0\leq
i\leq3$.
\end{Corollary}
\begin{proof}
Immediate from Corollary~\ref{cor:zbasis}.
\end{proof}

\section{Discussion}\label{Section16}

In this section we compare our main results with the results of Koornwinder~\cite{Koo1,Koo2}.

Recall from Def\/inition~\ref{def:udaha1} that $\hat H_q$ is the universal DAHA of type
$(C^\vee_1,C_1)$.
In~\cite{Koo1,Koo2} Koornwinder works with a~related algebra $\tilde{\mathfrak{H}}$ called the DAHA of type
$(C^\vee_1,C_1)$.
We will compare these algebras shortly.
Recall the set $\mathbb I=\lbrace0,1,2,3\rbrace$.
By Lemma~\ref{prop:altpres} the $\mathbb F$-algebra $\hat H_q$ has a~presentation by generators $\lbrace
t_i\rbrace_{i\in\mathbb I},\lbrace T_i\rbrace_{i\in\mathbb I}$ and relations
\begin{gather*}
t^2_i=T_i t_i-1,
\qquad
i\in\mathbb I,
\\
T_i
\ \text{is central}, \qquad
i\in\mathbb I,
\\
t_0t_1t_2t_3=q^{-1}.
\end{gather*}

\begin{Definition}
Let $\lbrace P_i\rbrace_{i\in\mathbb I}$ denote scalars in $\mathbb F$.
Def\/ine an $\mathbb F$-algebra $\hat H_q(P_0,P_1,P_2,P_3)$ by generators $\lbrace t_i\rbrace_{i\in\mathbb
I}$ and relations
\begin{gather*}
t^2_i=P_i t_i-1,
\qquad
i\in\mathbb I,
\qquad
t_0t_1t_2t_3=q^{-1}.
\end{gather*}
\end{Definition}

\begin{Lemma}
For $i\in\mathbb I$ the element $t_i$ of $\hat H_q(P_0,P_1,P_2,P_3)$ is invertible and $t_i+t^{-1}_i=P_i$.
\end{Lemma}

By construction there exists a~unique $\mathbb F$-algebra homomorphism $\hat H_q\to\hat
H_q(P_0,P_1,P_2,P_3)$ that sends $t_i\mapsto t_i$ and $T_i\mapsto P_i$ for $i\in\mathbb I$.
This map is surjective.
We denote this map by $\varepsilon(P_0,P_1,P_2,P_3)$.

Recall the elements $A$, $B$, $C$ of $\hat H_q$:
\begin{gather*}
A=t_1t_0+(t_1t_0)^{-1}=t_0t_1+(t_0t_1)^{-1},
\qquad
B=t_3t_0+(t_3t_0)^{-1}=t_0t_3+(t_0t_3)^{-1},
\\
C=t_2t_0+(t_2t_0)^{-1}=t_0t_2+(t_0t_2)^{-1}.
\end{gather*}
We retain the notation $A$, $B$, $C$ for their images under $\varepsilon(P_0,P_1,P_2,P_3)$.
Recall from Def\/ini\-tion~\ref{def:calx} the subalgebra $\mathbb T=\langle t^{\pm1}_0,T_1,T_2,T_3\rangle$ of
$\hat H_q$.
By~\eqref{eq:Ti}, $\mathbb T=\langle t_0,T_0,T_1,T_2,T_3\rangle$.
The subalgebra $\langle A,B,C,\mathbb T\rangle$ of $\hat H_q$ was discussed in
Propositions~\ref{prop:close},~\ref{prop:closer} and Theorems~\ref{thm:com},~\ref{thm:t0com}.

\begin{Definition}
Consider the subalgebra $\langle A,B,C,\mathbb T\rangle$ of $\hat H_q$.
Let $\mathcal A$ denote the image of $\langle A,B,C,\mathbb T\rangle$ under the map
$\varepsilon(P_0,P_1,P_2,P_3)$.
Thus $\mathcal A$ is the subalgebra of $\hat H_q(P_0,P_1,P_2,P_3)$ generated by $A$, $B$, $C$, $t_0$.
\end{Definition}
\begin{Proposition}\label{prop:myver}
The $\mathbb F$-algebra $\mathcal A$ is presented by generators and relations in the following way.
The generators are $A$, $B$, $C$, $t_0$.
The relations assert that $t_0$ is central and $t^2_0=P_0t_0-1$,
\begin{gather*}
A+\frac{qBC-q^{-1}CB}{q^2-q^{-2}}=\frac{\alpha}{q+q^{-1}},
\\
B+\frac{qCA-q^{-1}AC}{q^2-q^{-2}}=\frac{\beta}{q+q^{-1}},
\\
C+\frac{qAB-q^{-1}BA}{q^2-q^{-2}}=\frac{\gamma}{q+q^{-1}},
\\
q^{-1}ACB+q^{-2}A^2+q^{-2}B^2+q^2C^2-q^{-1}A\alpha-\;q^{-1}B\beta-qC\gamma
\\
\qquad
{}=
\big(q+q^{-1}\big)^2-\big(q^{-1}t_0+qt^{-1}_0\big)^2-P_1^2-P_2^2-P_3^2-\big(q^{-1}t_0+qt^{-1}_0\big)P_1P_2P_3,
\end{gather*}
where
\begin{gather*}
\alpha=\big(q^{-1}t_0+qt^{-1}_0\big)P_1+P_2P_3,
\qquad
\beta=\big(q^{-1}t_0+qt^{-1}_0\big)P_3+P_1P_2,
\\
\gamma=\big(q^{-1}t_0+qt^{-1}_0\big)P_2+P_3P_1,
\qquad
t^{-1}_0=P_0-t_0.
\end{gather*}
\end{Proposition}

\begin{proof}
In the relations of Theorem~\ref{thm:t0com}, f\/irst replace $t^{-1}_0$ by $T_0-t_0$ and then replace $T_i$
by $P_i$ for $i\in\mathbb I$.
\end{proof}

By the f\/irst three displayed relations in Proposition~\ref{prop:myver}, the $\mathbb F$-algebra $\mathcal
A$ is generated by $t_0$ together with any two of $A$, $B$, $C$.
We now give a~presentation of $\mathcal A$ by generators and relations, using the generators $A$, $B$,
$t_0$.
\begin{Proposition}
\label{prop:myverABt}
The $\mathbb F$-algebra $\mathcal A$ is presented by generators $A$, $B$, $t_0$ and relations
\begin{gather*}
t_0A=A t_0,
\qquad
t_0B=B t_0,
\qquad
t^2_0=P_0t_0-1,
\\
A^2B-\big(q^2+q^{-2}\big)ABA+BA^2+\big(q^2-q^{-2}\big)^2B+\big(q-q^{-1}\big)^2A\gamma
\\
\qquad
{}=\big(q-q^{-1}\big)\big(q^2-q^{-2}\big)\beta,
\\
B^2A-\big(q^2+q^{-2}\big)BAB+AB^2+\big(q^2-q^{-2}\big)^2A+\big(q-q^{-1}\big)^2B\gamma
\\
\qquad
{}=\big(q-q^{-1}\big)\big(q^2-q^{-2}\big)\alpha,
\\
q^{-1}A{\mathcal C}B+q^{-2}A^2+q^{-2}B^2+q^2{\mathcal C}^2-q^{-1}A\alpha-\;q^{-1}B\beta-q{\mathcal C}
\gamma
\\
\qquad
{}=\big(q+q^{-1}\big)^2-\big(q^{-1}t_0+qt^{-1}_0\big)^2-P_1^2-P_2^2-P_3^2-\big(q^{-1}t_0+qt^{-1}_0\big)P_1P_2P_3,
\end{gather*}
where
\begin{gather*}
\mathcal C=\frac{\gamma}{q+q^{-1}}-\frac{qAB-q^{-1}BA}{q^2-q^{-2}},
\qquad
\alpha=\big(q^{-1}t_0+qt^{-1}_0\big)P_1+P_2P_3,
\\
\beta=\big(q^{-1}t_0+qt^{-1}_0\big)P_3+P_1P_2,
\qquad
\gamma=\big(q^{-1}t_0+qt^{-1}_0\big)P_2+P_3P_1,
\qquad
t^{-1}_0=P_0-t_0.
\end{gather*}
\end{Proposition}

\begin{proof}
In the f\/irst two displayed relations of Proposition~\ref{prop:myver}, eliminate $C$ using the third
displayed relation.
\end{proof}

We now bring in the work of Koornwinder~\cite{Koo1, Koo2}.
In~\cite[equations~(3.1)--(3.4)]{Koo1} Koornwinder def\/ines an algebra $\tilde{\mathfrak H}$.
The def\/inition involves some scalars $q$, $a$, $b$, $c$, $d$.
For notational convenience we replace Koornwinder's $q$, $a$, $b$, $c$, $d$ by their squares.
\begin{Definition}[\protect{\cite[equations (3.1)--(3.4)]{Koo1}}] 
Fix nonzero scalars $a$, $b$, $c$, $d$ in $\mathbb F$.
The $\mathbb F$-algebra $\tilde{\mathfrak H}=\tilde{\mathfrak H}_q(a,b,c,d)$ is def\/ined by generators
$Z$, $Z^{-1}$, $\mathcal T_1$, $\mathcal T_0$ and relations
\begin{gather*}
\big(\mathcal T_1+a^2b^2\big)(\mathcal T_1+1)=0,
\qquad
\big(\mathcal T_0+q^{-2}c^2d^2\big)(\mathcal T_0+1)=0,
\\
\big(\mathcal T_1Z+a^2\big)\big(\mathcal T_1Z+b^2\big)=0,
\qquad
\big(q^2\mathcal T_0Z^{-1}+c^2\big)\big(q^2\mathcal T_0Z^{-1}+d^2\big)=0.
\end{gather*}
\end{Definition}

\begin{Lemma}[\protect{\cite[equations (3.5),~(3.6)]{Koo1}}] The elements $\mathcal T_1$, $\mathcal T_0$ of $\tilde{\mathfrak H}$
are invertible and
\begin{gather*}
\mathcal T^{-1}_1=-a^{-2}b^{-2}\mathcal T_1-1-a^{-2}b^{-2},
\qquad
\mathcal T^{-1}_0=-q^2c^{-2}d^{-2}\mathcal T_0-1-q^2c^{-2}d^{-2}.
\end{gather*}
\end{Lemma}

From now on, assume
\begin{gather}
P_0=ab+a^{-1}b^{-1},
\qquad
P_1=ab^{-1}+a^{-1}b,
\label{eq:P0P1}
\\
P_2=cd^{-1}+c^{-1}d,
\qquad
P_3=q^{-1}cd+qc^{-1}d^{-1}.
\label{eq:P2P3}
\end{gather}
\begin{Lemma}
\label{lem:terkor}
There exists an isomorphism of $\mathbb F$-algebras
\begin{gather*}
\hat H_q(P_0,P_1,P_2,P_3)\to\tilde{\mathfrak H}_q(a,b,c,d)
\end{gather*}
that sends
\begin{gather*}
t_0\mapsto-a^{-1}b^{-1}\mathcal T_1,
\qquad
t_1\mapsto-ab\mathcal T^{-1}_1Z^{-1},
\qquad
t_2\mapsto-q^{-2}cd Z\mathcal T^{-1}_0,
\qquad
t_3\mapsto-q c^{-1}d^{-1}\mathcal T_0.
\end{gather*}
The inverse isomorphism sends
\begin{gather*}
Z\mapsto qt_2t_3,
\qquad
Z^{-1}\mapsto t_0t_1,
\qquad
\mathcal T_1\mapsto-abt_0,
\qquad
\mathcal T_0\mapsto-q^{-1}cd t_3.
\end{gather*}
\end{Lemma}

\begin{proof}
One checks that the above maps are $\mathbb F$-algebra homomorphisms, and that they are inverses.
Consequently they are isomorphisms.
\end{proof}

From now on, we identify the $\mathbb F$-algebras $\hat H_q(P_0,P_1,P_2,P_3)$ and $\tilde{\mathfrak
H}_q(a,b,c,d)$ via the isomorphism in Lemma~\ref{lem:terkor}, and call the result $\tilde{\mathfrak H}$.

In~\cite[equations (3.8),~(3.9)]{Koo1} Koornwinder discusses two elements of $\tilde{\mathfrak H}$.
The f\/irst is $Y+q^{-2}a^2b^2c^2d^2Y^{-1}$ where $Y=\mathcal T_1\mathcal T_0$.
The second is $Z+Z^{-1}$.
These elements are related to $A$, $B$ as follows.
\begin{Lemma}
\label{lem:ZYvsAB}
In the algebra $\tilde{\mathfrak H}$,
\begin{gather*}
Z+Z^{-1}=A,
\qquad
Y+q^{-2}a^2b^2c^2d^2Y^{-1}=q^{-1}abcd B.
\end{gather*}
\end{Lemma}
\begin{proof}
Use Lemma~\ref{lem:terkor}.
\end{proof}

In~\cite[Def\/inition~6.1]{Koo1} Koornwinder def\/ines an $\mathbb F$-algebra $\widetilde{AW}(3,Q_0)$ by
generators and relations.
See also~\cite[Def\/inition~2.5]{Koo2}.
In~\cite[Corollary~6.3]{Koo1} Koornwinder displays an injection of $\mathbb F$-algebras
$\widetilde{AW}(3,Q_0)\to\tilde{\mathfrak H}$.
Consider the image of $\widetilde{AW}(3,Q_0)$ under this injection.
By construction and Lemma~\ref{lem:ZYvsAB}, this image is the subalgebra of $\tilde{\mathfrak H}$ generated
by $A$, $B$, $t_0$.
In other words, the image is $\mathcal A$.
Thus~\cite[Def\/inition~6.1, Corollary~6.3]{Koo1} yields a~presentation of~$\mathcal A$ by generators and
relations, using the generators $A$, $B$, $t_0$.
The presentation looks as follows in terms of $\lbrace P_i\rbrace_{i\in\mathbb I}$.
\begin{Theorem}
[\protect{\cite[Def\/inition~6.1, Corollary~6.3]{Koo1}}]
\label{prop:kverABt}
The $\mathbb F$-algebra $\mathcal A$ is presented by generators $A$, $B$, $t_0$ and relations
\begin{gather*}
t_0A=A t_0,
\qquad
t_0B=B t_0,
\qquad
t^2_0=P_0t_0-1,
\\
A^2B-\big(q^2+q^{-2}\big)ABA+BA^2+\big(q^2-q^{-2}\big)^2B+\big(q-q^{-1}\big)^2A\gamma
=\big(q-q^{-1}\big)\big(q^2-q^{-2}\big)\beta,
\\
B^2A-\big(q^2+q^{-2}\big)BAB+AB^2+\big(q^2-q^{-2}\big)^2A+\big(q-q^{-1}\big)^2B\gamma
=\big(q-q^{-1}\big)\big(q^2-q^{-2}\big)\alpha,
\\
\frac{ABAB}{(q^2-q^{-2})^2}-\frac{BABA(q^4+1+q^{-4})}{(q^2-q^{-2})^2}+\frac{B^2A^2\big(q^2+q^{-2}\big)}
{(q^2-q^{-2})^2}+A^2\big(q^2+q^{-2}\big)+B^2\big(q^2+q^{-2}\big)
\\
\qquad{}
+\frac{AB\gamma}{(q+q^{-1})^2}+\frac{BA\gamma\big(q\!-\!q^{-1}\big)\big(q^3\!-\!q^{-3}\big)}{(q^2\!-\!q^{-2})^2}\!
-\!\frac{A\alpha\big(q^3\!-\!q^{-3}\big)}{q^2\!-\!q^{-2}}\!
-\!\frac{B\beta\big(q^3\!-\!q^{-3}\big)}{q^2\!-\!q^{-2}}\!-\!\frac{\gamma^2}{(q+q^{-1})^2}
\\
\qquad{}
=\big(q+q^{-1}\big)^2-\big(q^{-1}t_0+qt^{-1}_0\big)^2-P_1^2-P_2^2-P_3^2-\big(q^{-1}t_0+qt^{-1}_0\big)P_1P_2P_3,
\end{gather*}
where
\begin{gather*}
\alpha=\big(q^{-1}t_0+qt^{-1}_0\big)P_1+P_2P_3,
\qquad
\beta=\big(q^{-1}t_0+qt^{-1}_0\big)P_3+P_1P_2,
\\
\gamma=\big(q^{-1}t_0+qt^{-1}_0\big)P_2+P_3P_1,
\qquad
t^{-1}_0=P_0-t_0.
\end{gather*}
\end{Theorem}

\begin{proof}
Write~\cite[Def\/inition~6.1]{Koo1} and~\cite[Corollary~6.3]{Koo1} in terms of $A$, $B$, $t_0$ and $\lbrace
P_i\rbrace_{i\in\mathbb I}$, using Lemma~\ref{lem:ZYvsAB} together with~\eqref{eq:P0P1},~\eqref{eq:P2P3}.
\end{proof}

In this paper we presented several subalgebras of $\hat H_q$ and $\tilde{\mathfrak H}$ by generators and
relations.
We now compare these presentations.
Theorems~\ref{thm:main1} and~\ref{thm:main3} together give a~presentation of the subalgebra $\langle
A,B,C\rangle$ of $\hat H_q$ by generators and relations, using the generators $A$, $B$, $C$.
Theorem~\ref{thm:t0com} gives a~presentation of the subalgebra $\langle A,B,C,\mathbb T\rangle$ of $\hat
H_q$ by generators and relations, using the generators $A$, $B$, $C$, $t^{\pm1}_0$, $T_1$, $T_2$, $T_3$.
Proposition~\ref{prop:myver} gives a~presentation of the subalgebra $\mathcal A$ of $\tilde{\mathfrak H}$
by generators and relations, using the generators $A$, $B$, $C$, $t_0$.
Proposition~\ref{prop:myverABt} and Theorem~\ref{prop:kverABt} each give a~presentation of $\mathcal A$ by
generators and relations, using the gene\-ra\-tors~$A$,~$B$,~$t_0$.
We now discuss the logical implications between Proposition~\ref{prop:myver},
Proposition~\ref{prop:myverABt}, and Theorem~\ref{prop:kverABt}.
Proposition~\ref{prop:myverABt} is discovered from Proposition~\ref{prop:myver} by partially
eliminating~$C$.
Proposition~\ref{prop:myver} is discovered from Proposition~\ref{prop:myverABt} and the knowledge that
$\mathcal C=C$.
Theorem~\ref{prop:kverABt} is discovered from Proposition~\ref{prop:myverABt} by eliminating~$\mathcal C$.
Proposition~\ref{prop:myverABt} is readily verif\/ied using Theorem~\ref{prop:kverABt}.
However Proposition~\ref{prop:myverABt} is not readily discovered using Theorem~\ref{prop:kverABt} alone.
Proposition~\ref{prop:myverABt} is discovered using Theorem~\ref{prop:kverABt} and the knowledge that
$\mathcal C$ simplif\/ies things.
Neither~$C$ nor $\mathcal C$ appears in~\cite{Koo1,Koo2}.

In~\cite{Koo2} Koornwinder discusses an algebra $S(\tilde{\mathfrak H})$ known as the spherical subalgebra
of $\tilde{\mathfrak H}$.
In~\cite[Theorem~3.2]{Koo2} Koornwinder displays an $\mathbb F$-algebra isomorphism ${\rm AW}(3,Q_0)\to
S(\tilde{\mathfrak H})$, where ${\rm AW}(3,Q_0)$ is the homomorphic image of $\widetilde{AW}(3,Q_0)$ described
in~\cite[Section~2]{Koo1}.
By~\cite[Section~3]{Koo2} the multiplicative identity of $S(\tilde{\mathfrak H})$ is a~certain
idempotent $P_{\rm sym}$ in $\tilde{\mathfrak H}$.
But $P_{\rm sym}\not=1$, so~$S(\tilde{\mathfrak H})$ and $\tilde{\mathfrak H}$ do not share the same~1.
Therefore $S(\tilde{\mathfrak H})$ is not a~subalgebra of $\tilde{\mathfrak H}$ according to our convention
from Section~\ref{Section1}.
As far as we know, the results of the present paper are unrelated to~$S(\tilde{\mathfrak H})$.

\subsection*{Acknowledgements}

The author thanks Kazumasa Nomura for giving this paper a~close reading and
of\/fering many valuable suggestions.

\pdfbookmark[1]{References}{ref}
\LastPageEnding

\end{document}